\documentclass[12pt,reqno]{amsart}

\usepackage{mathabx}
\usepackage{amssymb}
\usepackage{bbm}
\usepackage{mathrsfs}
\usepackage[all]{xy}
\usepackage{bm}
\usepackage{float}

\newtheorem{fac}{Fact}[section]
\newtheorem{theo}[fac]{Theorem}
\newtheorem{prop}[fac]{Proposition}
\newtheorem{lem}[fac]{Lemma}
\newtheorem{sub}[fac]{Sublemma}
\newtheorem{coro}[fac]{Corollary}

\theoremstyle{definition}
\newtheorem{ttt}[fac]{}
\newtheorem{defi}[fac]{Definition}

\theoremstyle{remark}
\newtheorem{ex}[fac]{Example}
\newtheorem{rem}[fac]{Remark}
\newtheorem{rems}[fac]{Remarks}

\newtheorem*{con}{Convention}
\newtheorem*{conv}{Conventions and notations}
\newtheorem*{org}{Organisation of the article}
\newtheorem*{qu}{Question}

\newcommand{\br}{ }
\newcommand{\brr}{, }

\newcommand{\Spec}{\mathop{\text{\rm Spec}}\nolimits}
\newcommand{\Gal}{\mathop{\text{\rm Gal}}\nolimits}
\newcommand{\Div}{\mathop{\text{\rm Div}}\nolimits}
\newcommand{\Pic}{\mathop{\text{\rm Pic}}\nolimits}
\newcommand{\Aut}{\mathop{\text{\rm Aut}}\nolimits}
\newcommand{\rk}{\mathop{\text{\rm rk}}\nolimits}
\renewcommand{\div}{\mathop{\text{\rm div}}\nolimits}
\newcommand{\pro}{\mathop{\text{\rm pr}}\nolimits}
\newcommand{\res}{\mathop{\text{\rm res}}\nolimits}
\newcommand{\im}{\mathop{\text{\rm im}}\nolimits}

\newcommand{\Br}{\mathop{\text{\rm Br}}\nolimits}
\newcommand{\Sym}{\mathop{\text{\rm Sym}}\nolimits}
\newcommand{\Sp}{\mathop{\text{\rm Sp}}\nolimits}
\newcommand{\Stab}{\mathop{\text{\rm Stab}}\nolimits}
\newcommand{\Split}{\mathop{\text{\rm Split}}\nolimits}
\newcommand{\Pol}{\mathop{\text{\rm Pol}}\nolimits}

\newcommand{\pr}{\text{\rm pr}}
\newcommand{\sep}{\text{\rm sep}}
\newcommand{\id}{\text{\rm id}}
\renewcommand{\l}{\text{\rm l}}
\newcommand{\ds}{\text{\rm ds}}
\newcommand{\maxi}{\text{\rm max}}
\newcommand{\sing}{\text{\rm sing}}

\newcommand{\bbF}{{\mathbbm F}}
\newcommand{\bbQ}{{\mathbbm Q}}
\newcommand{\bbZ}{{\mathbbm Z}}

\newcommand{\calD}{{\mathscr{D}}}
\newcommand{\calK}{{\mathscr{K}}}
\newcommand{\calL}{{\mathscr{L}}}
\newcommand{\calO}{{\mathscr{O}}}

\newcommand{\et}{{{\text{\rm {\'e}t}}}}

\newcommand{\Ab}{{\text{\bf A}}}
\newcommand{\Hb}{{\text{\bf H}}}
\newcommand{\Pb}{{\text{\bf P}}}

\newcounter{abc}
\newenvironment{abc}{\begin{list}{\rm \alph{abc}) }%
{\usecounter{abc} \leftmargin=0.0pt \labelsep=0.0pt %
\listparindent=0.0pt \labelwidth=0.0pt \parsep=\smallskipamount%
 \itemsep=0.0pt \topsep=0.0pt \partopsep=\smallskipamount}}{\end{list}}

\newcounter{iii}
\newenvironment{iii}{\begin{list}{\rm \roman{iii}) }%
{\usecounter{iii} \leftmargin=0.0pt \labelsep=0.0pt %
\listparindent=0.0pt \labelwidth=0.0pt \parsep=\smallskipamount%
 \itemsep=0.0pt \topsep=0.0pt \partopsep=\smallskipamount}}{\end{list}}

\makeatletter
\def\hsmash{\relax 
  \ifmmode\def\next{\mathpalette\mathhsm@sh}\else\let\next\makehsm@sh
  \fi\next}
\def\makehsm@sh#1{\setbox\z@\hbox{#1}\finhsm@sh}
\def\mathhsm@sh#1#2{\setbox\z@\hbox{$\m@th#1{#2}$}\finhsm@sh}
\def\finhsm@sh{\wd\z@\z@ \box\z@}
\makeatother

\def\rightend#1#2{{%
 \leavevmode\nobreak\hskip .5em plus 1fil
 \penalty600 \hskip 0pt plus -1filll
 \vadjust{}\nobreak\hskip 0pt plus 1filll%
 #1\parfillskip=#2\relax \par}}

\def\eop{\ifmmode\rule[-22pt]{0pt}{1pt}\ifinner\tag*{$\square$}\else\eqno{\square}\fi\else\rightend{$\square$}{0pt}\fi}

\thanks{}

\title[On plane quartics with a Galois invariant Steiner hexad]{On plane quartics with a Galois invariant \\Steiner hexad}

\begin{document}

\author{Andreas-Stephan Elsenhans}

\address{Institut f\"ur Mathematik\\ Universit\"at W\"urzburg\\ Emil-Fischer-Stra\ss e 30\\ D-97074 W\"urzburg\\ Germany}
\email{stephan.elsenhans@mathematik.uni-wuerzburg.de}
\urladdr{https://math.uni-paderborn.de/ag/ca/elsenhans/}

\author{J\"org Jahnel}

\address{\mbox{Department Mathematik\\ Univ.\ \!Siegen\\ \!Walter-Flex-Str.\ \!3\\ D-57068 \!Siegen\\ \!Germany}}
\email{jahnel@mathematik.uni-siegen.de}
\urladdr{http://www.uni-math.gwdg.de/jahnel}


\date{November~26,~2018}



\begin{abstract}
We describe a construction of plane quartics with prescribed Galois operation on the 28~bitangents, in the particular case of a Galois invariant Steiner~hexad. As~an application, we solve the inverse Galois problem for degree two del Pezzo surfaces in the corresponding particular case.
\end{abstract}

\maketitle
\thispagestyle{empty}

\section{Introduction}

It is well known~\cite{Pl} that a nonsingular plane quartic
curve~$C$
over an algebraically closed field of
characteristic~$\neq \!2$
has exactly 28 bitangents. The~same is still true if the base field is only separably closed, as is easily deduced from~\cite[Theorem~1.6]{Va}.
If~$C$
is defined over a separably non-closed
field~$k$
then the bitangents are defined over a finite extension field
$l$
of~$k$,
which is normal and separable, and permuted by the Galois
group~$\Gal(l/k)$.

By~far not every permutation may~occur. In~order to illustrate this, let us write
$\Omega_C$
for the set of all bitangents
of~$\smash{C_{\overline{k}}}$.
First~of all, every pair
$\smash{\{L, L'\} \in \Omega_C^{(2)}}$
of bitangents defines in a natural way a divisor~class
$$\Pi(\{L,L'\}) \in \Pic(C_{k^\sep})_2 = \Pic(C_{\overline{k}})_2 \cong H^1_\et(C_{\overline{k}}, \bbZ/2\bbZ) \, .$$
The~mapping
$\Pi\colon \Omega_C^{(2)} \to \Pic(C_{\overline{k}})_2$
is exactly six-to-one
onto~$\Pic(C_{\overline{k}})_2 \setminus \{\calO_{_{\overline{k}}}\}$,
the preimage of an element being classically called a Steiner~hexad. This~shows that a permutation
$\smash{\sigma \in \Sym \Omega_C \cong S_{28}}$
can be admissible only if the induced operation
$\smash{\sigma^{(2)} \in \Sym \Omega_C^{(2)} \cong S_{28\cdot27/2} = S_{378}}$
on 2-sets keeps the 63 Steiner hexads as a block~system.

Moreover,~$\Pic(C_{\overline{k}})_2$
is canonically equipped with the Weil pairing
$$\langle.\,,.\rangle\colon \Pic(C_{\overline{k}})_2 \times \Pic(C_{\overline{k}})_2 \longrightarrow \mu_2 \, .$$
An~admissible permutation therefore must provide an automorphism
of~$\Pic(C_{\overline{k}})_2$
that is symplectic with respect
to~$\langle.\,,.\rangle$.

It~turns out that the subgroups
$G_C \subset \Sym \Omega_C$
of all permutations that are admissible in the sense described are, independently of the choice
of~$C$,
permutation isomorphic to one and the same subgroup
$G \subset S_{28}$.
A~natural question arising is thus the~following.

\begin{qu}
Given~a
field~$k$
and a subgroup
$g \subseteq G$,
does there exist a nonsingular plane
quartic~$C$
over~$k$,
for which the group homomorphism
$\Gal(k^\sep/k) \to G \subset S_{28}$,
given by the Galois operation on the 28~bitangents, has the
subgroup~$g$
as its~image?
\end{qu}

\begin{rem}
One~has that
$G \cong \Sp_6(\bbF_{\!2})$.
This~is the simple group of order
$1\,451\,520$.
Up~to conjugation, it has only one subgroup of index 28, which is a maximal~subgroup. Hence,~the permutation group
$G \subset S_{28}$
is~primitive. It~is even
\mbox{$2$-transitive}
in view of \cite[Theorem~7.7A]{DM}.

On~the other hand, as a permutation group,
$G \subset S_{28}$
is self-normalising. Indeed, the centraliser
$C_{S_{28}}(G)$
is trivial by \cite[Theorem~4.2A.(vi)]{DM} and
$\Sp_6(\bbF_{\!2})$
is known to have only inner automorphisms.
This shows that the permutation isomorphisms
$\smash{G_C \stackrel{\cong}{\rightarrow} G}$
are uniquely determined up to conjugation by elements
of~$G$.
In~particular, the question above depends only on the conjugacy class of the
subgroup~$g \subseteq G$.
\end{rem}

The~group
$G \cong \Sp_6(\bbF_{\!2})$
has
$1369$~con\-ju\-gacy classes of subgroups. Among~these, there are eight maximal subgroups, which are of indices
$28$,
$36$,
$63$,
$120$,
$135$,
$315$,
$336$,
and
$960$,~respectively.

\looseness-1
An~example of a nonsingular plane quartic
over~$\bbQ$
such that
$\Gal(\overline\bbQ/\bbQ)$
operates on the bitangents via the full
$\Sp_6(\bbF_{\!2})$
has been constructed by T.~Shioda~\cite[Section~3]{Shi} in 1993 and, almost at the same time, by R.~Ern\'e~\cite[Corollary~3]{Er}. Moreover,~there is an obvious approach to construct examples for the groups contained in the index
$28$~subgroup.
Indeed,~in this case, there is a rational~bitangent. One~may start with a cubic surface with the right Galois operation~\cite{EJ15}, blow-up a rational point, and use the connection between degree two del Pezzo surfaces and plane quartics~\cite[Theorem~3.3.5]{Ko96}, cf.\ the application discussed~below.\smallskip

In~this article, we deal with the subgroup
$U_{63} \subset G$
of
index~$63$
and the groups contained~within. More precisely, we show the following result, which answers, in the case that the base field is infinite and of characteristic
$\neq \!2$,
a more refined question than the one asked~above.

\begin{theo}
\label{general}
Let an infinite
field\/~$k$
of characteristic
not\/~$2$,
a normal and separable extension
field\/~$l$,
and an injective group~homomorphism
$$i\colon \Gal(l/k) \hookrightarrow U_{63}$$
be~given. Then there exists a nonsingular quartic
curve\/~$C$
over\/~$k$
such that\/
$l$
is the field of definition of the 28 bitangents and each\/
$\sigma \in \Gal(l/k)$
permutes the bitangents as described by\/
$i(\sigma) \in G \subset S_{28}$.
\end{theo}

Among the 1369 conjugacy classes of subgroups
of~$G \cong \Sp_6(\bbF_{\!2})$,
1155 are contained
in~$U_{63}$.
By~which we mean that they have a member that is contained in a fixed subgroup of
index~$63$.

The~maximal subgroup
$U_{63}$
has a geometric meaning. 
Namely,~$G \cong \Sp_6(\bbF_{\!2})$
operates transitively on the 63 elements of
$\smash{\Pic(C_{\overline{k}})_2 \setminus \{\calO_{C_{\overline{k}}}\}}$.
Thus,~the subgroup
$U_{63} \subset G$
is just the point~stabiliser and the inclusion
$\Gal(l/k) \hookrightarrow U_{63}$
expresses the fact that there is a
\mbox{$k$-rational}
divisor class
in~$\Pic(C)_2$.
Consequently,~there is a Galois invariant Steiner hexad
on~$C$,~too.

Furthermore,~in this case, one of the corresponding del Pezzo surfaces of degree two has a
\mbox{$k$-rational}
conic~bundle. In~fact, there are two such bundles, which are mapped to each other under the Geiser~involution. We~use conic bundles for our proof of existence, which is completely~constructive.

\begin{rems}
\label{mark}
\begin{iii}
\item
As~noticed above, a
quartic~$C$
that provides a solution for a homomorphism
$i\colon \Gal(l/k) \hookrightarrow G$
may, as well, serve as a solution for any homomorphism
$\varphi_g \!\circ\! i$
differing
from~$i$
by an inner automorphism
$\varphi_g$
of~$G$,
for
some~$g \in G$.

In~Theorem~\ref{general}, we require
$\im i$
to be contained in the self-normalising subgroup
$U_{63} \subset G$.
Thus,~in general, one may expect that
$\varphi_g \!\circ\! i$
maps
to~$U_{63}$
only
for~$g \in U_{63}$.
I.e.,~we may disturb
$i$
by the inner automorphisms
of~$U_{63}$.
\item
However,~our approach is based on conic bundles with six split~fibres. It~provides, a priori, a quartic with the right Galois operation, up to an inner automorphism
of~$S_2 \wr S_6$.
As,~however,
$U_{63} \cong (S_2 \wr S_6) \!\cap\! A_{12} \subsetneqq S_2 \wr S_6$,
a rather sophisticated monodromy argument is necessary in order to complete the proof for the main~result.
\end{iii}
\end{rems}

As~an application, one may answer the analogous question for degree two del Pezzo surfaces. The~double cover
of~$\Pb^2$,
ramified at a nonsingular quartic
curve~$C$,
is a del Pezzo surface of degree two.
Here,~considerations can be made that are very similar to the ones~above. First~of all, it is well known~\cite[Theorem~26.2.(iii)]{Ma} that a del Pezzo
surface~$S$
of degree two over an algebraically closed field contains exactly 56 exceptional curves, i.e.\ such of self-intersection
number~$(-1)$.
Again,~the same is true when the base field is only separably~closed.
If~$S$
is defined over a separably non-closed
field~$k$
then the exceptional curves are defined over a normal and separable finite extension field
$l$
of~$k$
and permuted
by~$\Gal(l/k)$.
Once~again, not every permutation may~occur. The~maximal subgroup
$\smash{\widetilde{G} \subset \Sym \widetilde\Omega_C \cong S_{56}}$
that respects the intersection pairing is isomorphic to the Weil
group~$W(E_7)$~\cite[Theorem~23.9]{Ma}.

Every~bitangent
of~$C$
is covered by exactly two of the exceptional curves
of~$S$.
Thus,~for the operation
of~$\smash{\Gal(k^\sep/k)}$
on the 56 exceptional curves
on~$S$,
there seem to be two independent conditions. On~one hand,
$\smash{\Gal(k^\sep/k)}$
must operate via a subgroup
of~$\smash{W(E_7) \cong \widetilde{G} \subset \Sym \widetilde\Omega_C}$.
On~the other hand, the induced operation on the blocks of size two must take place via a subgroup of
$\Sp_6(\bbF_{\!2}) \cong G \subset \Sym \Omega_C \cong S_{28}$.
It~turns out, however, that there is an isomorphism
$\smash{W(E_7)/Z \stackrel{\cong}{\longrightarrow} \Sp_6(\bbF_{\!2})}$,
for
$Z \subset W(E_7)$
the centre, that makes the two conditions~equivalent.

The~group
$\smash{\widetilde{G} \cong W(E_7)}$
already has
$8074$
conjugacy classes of~subgroups. Two~subgroups with the same image under the quotient~map
$\smash{p\colon \widetilde{G} \twoheadrightarrow \widetilde{G}/Z \stackrel{\cong}{\longrightarrow} G}$
correspond to del Pezzo surfaces of degree two that are quadratic twists of each~other. Theorem~\ref{general} therefore extends word-by-word to del Pezzo surfaces of degree two and homomorphisms
$\smash{\Gal(l/k) \hookrightarrow \widetilde{G}}$
with image contained
in~$p^{-1}(U_{63})$.\smallskip

There~are further applications of Theorem~\ref{general} that concern cubic~surfaces. For~instance, we refine our previous result~\cite{EJ10} on the existence of cubic surfaces with a Galois invariant double six and generalise it from
$\bbQ$
to an arbitrary infinite field of characteristic
not~$2$.

\begin{org}
Section~\ref{zwei} summarises several general results on plane quartics, degree
$2$
del Pezzo surfaces, Steiner hexads, and conic bundles, which are necessary for our~arguments. They~are without doubt well-known to experts. Thus,~concerning this part, we do not claim any originality, except for the~presentation. Section~\ref{drei}~then describes our approach to the construction of plane quartics in detail, thereby proving Theorem~\ref{general}. The~explicit description of a conic bundle with six split fibres that are acted upon by Galois in a prescribed way, given in Proposition~\ref{cb_gal}, is the heart of this~approach. We~deduce the main result in Section~\ref{vier} by providing the necessary monodromy argument, which is related to the fact that
$U_{63} \cong (S_2 \wr S_6) \!\cap\! A_{12} \subsetneqq S_2 \wr S_6$.
Degree~two del Pezzo surfaces and their quadratic twists of are discussed afterwards, in Section~\ref{fuenf}, and, finally, we present the applications to cubic surfaces of particular types in Sections~\ref{sechs} and~\ref{sieben}. All calculations are with {\tt magma}~\cite{BCP}.
\end{org}

\begin{rems}
\begin{iii}
\item
Theorem~\ref{general} is clearly not true, in general, when
$k$
is a finite~field. For~example, there cannot be a nonsingular quartic curve
over~$\bbF_{\!3}$,
all whose bitangents are
\mbox{$\bbF_{\!3}$-rational},
simply because the projective plane contains only 13 
\mbox{$\bbF_{\!3}$-rational}
lines.
\item
We~ignore about
characteristic~$2$
in this article, as this case happens to be very different. Even~over an algebraically closed field, a plane quartic cannot have more than seven bitangents~\cite[p.~60]{SV}.
\end{iii}
\end{rems}

\begin{conv}
In~this article, we follow standard conventions and notations from Algebra and Algebraic Geometry, except for the~following.

\begin{iii}
\item
By a {\em field,} we mean a field of
characteristic~$\neq\! 2$.
For~the convenience of the reader, the assumption on the characteristic will be repeated in the formulations of our final results, but not during the intermediate~steps.
\item
When
$V$
is a finite-dimensional vector space over a
field~$k$,
then we denote its associated affine space
$\Spec \Sym V^\vee$
in the category of schemes
by~$\Ab(V)$.
The~elements
of~$V$
are then naturally in bijection with the
\mbox{$k$-rational}
points
on~$\Ab(V)$.

When~$E$
is an affine space, acted upon transitively and freely by the vector
space~$V$,
we put
$\Ab(E) := \Spec \Pol(E)$,
for
$\Pol(E)$
the algebra of polynomial functions
on~$E$.
This~is simply the pull-back of the
\mbox{$k$-algebra}
$\Sym V^\vee$
under any of the \mbox{bijections}
$\iota_{v_0}\colon E \to V, \,v_0+v \mapsto v$.
\item
For~a field
$k$
and an integer
$d>0$,
we write
$k[T]_d$
for the set of all monic polynomials of
degree~$d$
with coefficients
in~$k$.
This~is an affine space under the
\mbox{$d$-d}i\-men\-sional
\mbox{$k$-vector}
space of the polynomials of
degree~$< \!d$.

\end{iii}
\end{conv}

\section{Generalities on plane quartics and degree two \\del Pezzo surfaces}
\label{zwei}

\begin{defi}
Let
$\smash{C \subset \Pb^2_k}$
be a plane curve over a
field~$k$.

\begin{iii}
\item
Then,~by a {\em contact conic\/}
of~$C$,
one means a conic
$D \subset \Pb^2_k$
such that, at every geometric point
$\smash{p \in (C \cap D)(\overline{k})}$
of intersection, the intersection multiplicity is~even.
\item
With a contact
conic~$D$,
one associates an invertible sheaf
$\calD \in \Pic(C)_2$,
i.e.~one that is annihilated
by~$2$
in the Picard group. Just~put
$$\textstyle \calD := \calO(\frac12 C.D) \!\otimes\! \calO_{\Pb^2}(-1)|_C \, ,$$
for
$C.D$
the intersection divisor of
$D$
with~$C$.
\end{iii}
\end{defi}

Let~$C$
be a plane quartic over an algebraically closed field. Then~the genus
of~$C$
is
$\smash{g = \frac{3\cdot2}2 = 3}$,
which implies that
$\Pic(C)_2$
is a
$6$-dimensional vector space
over~$\bbF_{\!2}$.
In~this situation, the union
$L \cup L'$
of two bitangents
$L \neq L'$
is a degenerate contact~conic.
If~$L$
and~$\smash{L'}$
touch
$C$
in
$p_1$
and~$p_2$
and
$\smash{p'_1}$
and~$\smash{p'_2}$,
respectively, then the associated invertible sheaf
in~$\Pic(C)_2$~is
\begin{eqnarray*}
\textstyle \calO_C(\frac12 C.L) \!\otimes\! \calO_C(\frac12 C.L') \otimes \calO_{\Pb^2}(-1)|_C & = & \textstyle \calO_C(\frac12 C.L - \frac12 C.L') \\
 & = & \calO_C((p_1) + (p_2) - (p'_1) - (p'_2)) \\
 & = & \calO_C((p'_1) + (p'_2) - (p_1) - (p_2)) \, .
\end{eqnarray*}
This invertible sheaf is automatically nontrivial. Indeed, otherwise
$\calO_C((p_1) + (p_2))$
would have a non-constant section. However,~one has
$\calO_C(2(p_1) + 2(p_2)) \cong \calO_{\Pb^2}(1)|_C$
and the fact that
$H^1(\Pb^2, \calO_{\Pb^2}(-3)) = 0$
ensures that the latter sheaf has no global sections other than restrictions of global linear forms. Hence,
\begin{equation}
\label{lin_schnitt}
\textstyle \Gamma(C, \calO_C(2(p_1) + 2(p_2))) = \langle \frac{T_0}l, \frac{T_1}l, \frac{T_2}l \rangle \, ,
\end{equation}
for
$T_0$,
$T_1$,
and~$T_2$
the coordinate functions
on~$\Pb^2$
and
$l$
the linear form
defining~$L$.
Since,~as a rational function
on~$C$,
none of the linear combinations has only simple poles, one sees~that
\begin{equation}
\label{theta_odd}
h^0(C, \calO_C((p_1) + (p_2)))= 1 \, .
\end{equation}

\begin{rem}
\label{bita_div}
As~the conclusion of the considerations just made, one obtains a canonical~mapping
\begin{eqnarray*}
\Pi\colon \Omega_C^{(2)} = \{\text{pairs of bitangents of~} C\} & \longrightarrow & \Pic(C)_2 \setminus \{\calO_C\} \, , \\
\{L,L'\} & \mapsto & \textstyle \calO_C(\frac12 C.L - \frac12 C.L')
\end{eqnarray*}
from the
$\smash{\frac{28\cdot27}2 = 378}$
pairs of bitangents of
$C$
to the
$63$
nonzero elements
in~$\Pic(C)_2$.
\end{rem}

\begin{prop}
\label{detf}
Let\/~$C \subset \Pb^2_k$
be a nonsingular plane quartic over a
field\/~$k$.
Assume that there is given an invertible sheaf\/
$\calL \in \Pic(C)_2 \setminus \{\calO_C\}$,
i.e.\ one that is nontrivial,
\mbox{$2$-torsion}
in the Picard group, and defined
over\/~$k$.

\begin{abc}
\item
Then\/~$C$
may be written in the symmetric determinantal~form
$$C\colon q^2 - q_1 q_2 = 0 \, ,$$
for
$q$,
$q_1$,
and\/~$q_2$
three quadratic forms with coefficients
in\/~$k$.
\item
Furthermore, the equations\/
$q_1 = 0$
and\/~$q_2 = 0$
define contact conics
of\/~$C$
that are associated
with\/~$\calL$.
\item
All contact conics associated
with\/~$\calL$
form a one-dimensional family. They are of the type\/
$K_{(s:t)}\colon s^2 q_1 + 2st q + t^2 q_2 = 0$,
for\/
$(s:t) \in \Pb^1$.
\item
The double cover\/
$S\colon w^2 = q^2 - q_1 q_2$
carries two\/
\mbox{$k$-rational}
conic bundles, the projections of which down
to\/~$\Pb^2_k$
coincide with the
family\/~$K_{(s:t)}$.
\end{abc}\smallskip

\noindent
{\bf Proof.}
{\em
a) and~b)
Since~$C$
is a plane quartic and nonsingular, the adjunction formula shows that
$\smash{\calO_{\Pb^2}(1)|_C = \calK_C}$
is the canonical~sheaf. 
Clearly,~one
has~$\smash{\deg \calO_{\Pb^2}(1)|_C = 4}$.
Moreover,
$\calL^\vee$
is an invertible sheaf that is nontrivial and of
degree~$0$,
so that it has no non-zero~section.
Therefore, the Theorem of Riemann-Roch shows~that
$$h^0(C, \calL \!\otimes\! \calO_{\Pb^2}(1)|_C) = h^0(C, \calL\!\otimes\! \calK_C) - h^0(C, \calL^\vee) = 4 + 1 - g = 2 \, .$$
I.e., there is a pencil of effective divisors defining the invertible sheaf
$\smash{\calL \!\otimes\! \calO_{\Pb^2}(1)|_C}$.
Let
$(p_1) + \cdots + (p_4)$
be such an effective divisor for
$\smash{\calL \!\otimes\! \calO_{\Pb^2}(1)|_C}$
and
$(p_5) + \cdots + (p_8)$
be another. Then
$$2(p_1) + \cdots + 2(p_4) \, , \quad 2(p_5) + \cdots + 2(p_8) \, ,  \quad\text{and}\quad (p_1) + \cdots + (p_8)$$
are three effective divisors defining
$\smash{(\calL \!\otimes\! \calO_{\Pb^2}(1)|_C)^{\otimes 2} = \calO_{\Pb^2}(2)|_C}$.
Thus, one~has

\begin{iii}
\item
a contact conic
$q_1$
such that
$\div(q_1) = 2(p_1) + \cdots + 2(p_4)$,
\item
a contact conic
$q_2$
such that
$\div(q_2) = 2(p_5) + \cdots + 2(p_8)$,
and
\item
a conic
$q$
such that
$\div(q) = (p_1) + \cdots + (p_8)$.
\end{iii}

\noindent
Altogether,
$q^2, q_1 q_2 \in \Gamma(C, \calO_{\Pb^2}(4)|_C)$
both have
$2(p_1) + \cdots + 2(p_8)$
as its associated divisor. Therefore,~they must agree up to a constant
factor~$c \in k^*$.
In~other words,
$q^2 - cq_1q_2 = 0$
holds
on~$C$.
One~may normalise
$q_1$
so that this equation takes the~form
\begin{equation}
\label{eq1}
q^2 - q_1q_2 = 0 \, .
\end{equation}

In~order to make sure that this is the equation of the curve
$C$,
one still needs to exclude that (\ref{eq1}) holds identically on the whole
of~$\Pb^2$.
Assuming~it would hold identically, we first observe that
$q_1$
and~$q_2$
are not associates, since
$\div(q_1) \neq \div(q_2)$.
In~view of this, equation~(\ref{eq1}) implies that
$q$
splits into two non-associate linear
factors,~$q = l_1l_2$.
The~equation
$l_1^2 l_2^2 = q_1 q_2$
then shows, finally, that both
$q_1$
and~$q_2$
must define double lines, e.g.\
$q_1 \sim l_1^2$
and
$q_2 \sim l_2^2$.
But~now
$\div(q_1) = 2(p_1) + \cdots + 2(p_4)$
yields
$\div(l_1) = (p_1) + \cdots + (p_4)$,
so
$(p_1) + \cdots + (p_4)$
is a divisor defining
$\smash{\calO_{\Pb^2}(1)|_C}$
and not
$\smash{\calL \!\otimes\! \calO_{\Pb^2}(1)|_C}$,
a~contradiction.\smallskip

\noindent
c)
The conics
$\smash{D_{(s:t)} \colon s^2 q_1 + 2st q + t^2 q_2 = 0}$
are indeed contact conics
of~$C$.
For~$s=0$,
this was shown above. Otherwise,
on~$D_{(s:t)}$,
one has
$\smash{q_1 = -\frac{2t}s q - \frac{t^2}{s^2} q_2}$,
hence the equation
of~$C$
takes the form
$\smash{0 = q^2 - (-\frac{2t}s q - \frac{t^2}{s^2} q_2)q_2 = (q + \frac{t}s q_2)^2}$.
Moreover, the~contact conics
$\smash{D_{(s:t)}}$
are all associated with the same
$2$-torsion
invertible sheaf, as the base scheme
$\Pb^1$
is connected. According~to b), this sheaf is
exactly~$\calL$.

On~the other hand, we showed above that there is only a one-dimensional linear system of effective divisors that
define~$\smash{\calL \!\otimes\! \calO_{\Pb^2}(1)|_C}$.
Furthermore,~since
$2 < \deg C = 4$,
the divisor determines the conic uniquely, which implies the~claim.\smallskip

\noindent
d)
A direct calculation shows that the~conics
$$\textstyle w = \frac{s}t q_1 + q, \quad s^2 q_1 + 2st q + t^2 q_2 = 0$$
lie
on~$S$,
for~$(s\!:\!t) \in \Pb^1$.
For~$t=0$,
this is supposed to mean
$w=q$
and~$q_1 = 0$.
The~second conic bundle is obtained replacing
$w$
by~$-w$.
}
\eop
\end{prop}

\begin{rem}
Part~a) of this result is essentially \cite[Theorem~6.2.3]{Do}, where it is deduced from the general framework of theta characteristics. For~our purposes, this generality is neither necessary nor does it lead to additional~clarity. On~the other hand, we find the direct proof, just using the Riemann-Roch Theorem, quite instructive.
\end{rem}

\begin{coro}
Let\/~$C$
be a nonsingular plane quartic over an algebraically closed field. Then~the mapping\/
$\Pi\colon \{\text{pairs of bitangents of~} C\} \to \Pic(C)_2 \setminus \{\calO_C\}$
is surjective and precisely six-to-one.\medskip

\noindent
{\bf Proof.}
{\em
Let~$\calL \in \Pic(C)_2 \setminus \{\calO_C\}$
be any element. Then~the contact conics associated with
$\calL$
form a one-dimensional family, which can be written in the form
$\smash{D_{(s:t)} \colon s^2 q_1 + 2st q + t^2 q_2 = 0}$.
The quadratic forms
$q_1$,
$q_2$,
and~$q$
occurring may be described by symmetric
$3\times3$-matrices
$M_1$,
$M_2$,
and~$M$,
respectively. Degenerate~conics occur at the zeroes of
$\det(s^2 M_1 + 2st M + t^2 M_2)$,
which is a binary form that is homogeneous of degree~six. Consequently, not more than six of the conics may be degenerate and not more than six pairs of bitangents may be mapped
to~$\calL$
under~$\Pi$.
Since~$378 = 63\cdot6$,
the proof is~complete.
}
\eop
\end{coro}

\begin{defi}[cf.~{\cite[Section~6.1.1]{Do}}]
\label{St_hex}
Let~$\calL \in \Pic(C)_2 \setminus \{\calO_C\}$
be any element. Then~the preimage
$\Pi^{-1}(\calL)$
is called a {\em Steiner~hexad}.
\end{defi}

\begin{ttt}[Degree two del Pezzo surfaces]
Let
$C\colon Q=0$
be a nonsingular plane quartic over a
field~$k$.
Then,~for every
$\lambda \in k^*$,
there is the double cover
of~$\Pb^2$,
ramified
at~$C$,
given by
$S\colon \lambda w^2 = Q$.
This~is a del Pezzo surface of degree two~\cite[Theorem~III.3.5]{Ko96}. It~is equipped with the projection
$\pi\colon S \to \Pb^2$
and the Geiser~involution
$$g\colon (w;T_0\!:\!T_1\!:\!T_2) \mapsto (-w;T_0\!:\!T_1\!:\!T_2) \, .$$
Over each of the 28 bitangents
of~$C_{\overline{k}}$,
the double cover splits into two components. This yields exactly the 56 exceptional curves
on~$S_{\overline{k}}$.
In~fact,
if~$l = 0$
defines a bitangent then, modulo the linear
form~$l$,
the equation
of~$C$
is
$cq^2 = 0$,
for
$\smash{c \in \overline{k}^*}$
and a quadratic
form~$q$.
Thus,~the equation
of~$S$
takes the form
$\lambda w^2 = cq^2$,
which shows the splitting~into
\begin{equation}
\label{split_loc}
w = \pm\sqrt{\frac{c}\lambda} q \, .
\end{equation}
\end{ttt}

\begin{rem}[The blown-up model]
A del Pezzo surface of degree two over an algebraically closed field is isomorphic
to~$\Pb^2$,
blown up in seven points
$x_1, \ldots, x_7$
in general position~\cite[Theorem~24.4.(iii)]{Ma}. In~the blown-up model, the 56 exceptional curves are given as follows, cf.~\cite[Theorem~26.2]{Ma}.

\begin{iii}
\item
$e_i$,
for~$i=1, \ldots, 7$,
the inverse image of the blow-up point
$x_i$.
\item
$l_{ij}$,
for~$1 \leq i < j \leq 7$,
the strict transform of the line through
$x_i$
and~$x_j$.
The~class
of~$l_{ij}$
in
$\Pic(S)$
is~$l-e_i-e_j$.
\item
$\smash{\widetilde{l}_{ij}}$,
for~$1 \leq i < j \leq 7$,
the strict transform of the conic through all blow-up points excluding
$x_i$
and~$x_j$.
The~class
of~$\smash{\widetilde{l}_{ij}}$
in
$\Pic(S)$
is~$2l-e_1-\cdots-e_7 + e_i+e_j$.
\item
$\smash{\widetilde{e}_i}$,
for~$i=1, \ldots, 7$,
the strict transform of the unique singular cubic curve through all seven blow-up points that has
$x_i$
as a double~point. The~class
of~$\smash{\widetilde{e}_i}$
in
$\Pic(S)$
is
$3l-e_1-\cdots-e_7-e_i$.
\end{iii}

\noindent
The Geiser involution
$g$
interchanges
$e_i$
with
$\smash{\widetilde{e}_i}$,
for~$i=1, \ldots, 7$,
and
$\smash{l_{ij}}$
with
$\smash{\widetilde{l}_{ij}}$,
for
$1 \leq i < j \leq 7$.
Furthermore,~the maximal subgroup
$\smash{\widetilde{G} \subset S_{56}}$
that respects the intersection pairing operates transitively on systems of
$i$
mutually skew exceptional curves,
for~$i = 1, \ldots, 5$,
or~$7$~\cite[Corollary~26.8.(i)]{Ma}.
Thus,~given such a system of
$i$
exceptional curves, one may assume without restriction that the curves are
$e_1, \ldots, e_i$.
\end{rem}

\begin{ttt}[The Picard group]
\label{Picardgr}
The Picard group
$\Pic(S)$
of a degree two del Pezzo
surface~$S$
over an algebraically closed field is a free abelian group of
rank~$8$
and generated by the 56 exceptional curves. One~has the transfer map
$\pi_*\colon \Pic(S) \to \Pic(\Pb^2)$.
The~kernel
$\Pic(S)[\pi_*]$~is

\begin{iii}
\item
generated by all differences
$\calO_S(E - E')$
of exceptional curves,
\item
generated up to index three by
$\calO_S(E_2-E_1), \ldots, \calO_S(E_7-E_1)$,
and
$\smash{\calO_S(\widetilde{E}_1-E_1)}$,
for seven mutually skew exceptional curves
$E_1, \ldots, E_7$
and
$\smash{\widetilde{E}_1 := g(E_1)}$.
\end{iii}

\noindent
The second assertion easily follows from a short consideration in the blown-up~model.
\end{ttt}

\begin{ttt}[The various restriction homomorphisms]
For~the double cover
$S$
of~$\Pb^2$,
ramified at the nonsingular
quartic~$C$,
one also has the restriction homomorphism
$r\colon \Pic(S) \to \Pic(C)$.
Formula~(\ref{split_loc}) shows that, if
$E$
is an exceptional curve
on~$S$
and
$\pi(E)$
touches~$C$
in
$p_1$
and
$p_2$~then
$$\textstyle r(\calO_S(E)) = \calO_C(\frac12 C.\pi(E)) = \calO_C((p_1)+(p_2)) \, .$$
Consequently,
$\smash{r(\calO_S(2E)) = \calO_{\Pb^2}(1)|_C}$,
for every exceptional
curve~$E$.

Let~us denote the restriction of
$r$
to
$\Pic(S)[\pi_*]$
by
$r'\colon \Pic(S)[\pi_*] \to \Pic(C)$.
Then
\begin{equation}
\label{r_pi}
r'(\calO_S(E - E')) = \Pi(\{\pi(E), \pi(E')\}) \, .
\end{equation}
This~shows that
$r'$
maps
$\Pic(S)[\pi_*]$
only
to~$\Pic(C)_2$.
Moreover,~the homomorphism
$r'\colon \Pic(S)[\pi_*] \to \Pic(C)_2$
is obviously~onto. We have the induced surjective homomorphism
\begin{equation}
\label{restr}
\overline{r}\colon \Pic(S)[\pi_*]/2\Pic(S)[\pi_*] \longrightarrow \Pic(C)_2 \, .
\end{equation}

Finally,~recall that, via the Chern class homomorphism, there is a canonical isomorphism
$\Pic(S)/2\Pic(S) \cong H^2_\et(S, \bbZ/2\bbZ)$.
\end{ttt}

\begin{lem}[Criterion for Steiner hexads]
\label{crit_hex}
Let\/~$C \subset \Pb^2$
be a nonsingular plane~quartic over an algebraically closed field and\/
$S$
be the double cover
of\/~$\Pb^2$,
ramified~at\/~$C$.
Furthermore,~let\/
$\{E_1,E'_1,\ldots,E_6,E'_6\}$
be a set of twelve exceptional curves
on\/~$S$
such~that

\begin{iii}
\item
$E_1,\ldots,E_6$
are mutually skew,
\item
$E'_1,\ldots,E'_6$
are mutually skew, and
\item
$E_i \!\cdot\! E'_j = \delta_{ij}$,
for\/
$1 \leq i,j \leq 6$
and\/
$\delta_{ij}$
the Kronecker~symbol.
\end{iii}

\noindent
Then\/
$\{(\pi(E_1),\pi(E'_1)), \ldots, (\pi(E_6),\pi(E'_6))\}$
is a Steiner~hexad.\medskip

\noindent
{\bf Proof.}
{\em
The exceptional curves
$E_1,E'_1,\ldots,E_6,E'_6$
generate, together with the canonical class
$[K]$,
the whole Picard group
$\Pic(S)$
up to finite~index. Indeed,~the
$13\times13$
intersection~matrix
\begin{equation}
\label{rank8}
\left(
\begin{array}{rrrrrrrr}
-1 & 1 & 0 & 0 & \phantom{j}\ldots & 0 & 0 & -1 \\
 1 &-1 & 0 & 0 & \ldots & 0 & 0 & -1 \\
 0 & 0 &-1 & 1 & \ldots & 0 & 0 & -1 \\
 0 & 0 & 1 &-1 & \ldots & 0 & 0 & -1 \\
\\
 0 & 0 & 0 & 0 & \ldots &-1 & 1 & -1 \\
 0 & 0 & 0 & 0 & \ldots & 1 &-1 & -1 \\
-1 &-1 &-1 &-1 & \ldots &-1 &-1 &  2
\end{array}
\right)
\end{equation}
is of
rank~$8$.
(The index of
$\langle E_1,E'_1,\ldots,E_6,E'_6,[K] \rangle \subset \Pic(S)$
is,
in fact,~$2$.)

Thus,~we may establish the equality
$\calO_S(E_1+E'_1) = \cdots = \calO_S(E_6+E'_6)$
in
$\Pic(S)$
by just noting that each of the six divisors has intersection
number~$0$
with
$E_1,E'_1,\ldots,E_6,E'_6$
and~$(-2)$
with~$[K]$.
Consequently,
$$r(\calO_S(E_1+E'_1)) = \cdots = r(\calO_S(E_6+E'_6)) \, .$$
As~$\smash{r(\calO_S(2E)) = \calO_{\Pb^2}(1)|_C}$
for every exceptional curve, formula (\ref{r_pi}) yields~that
$$\Pi(\{\pi(E_1), \pi(E'_1)\}) = \cdots = \Pi(\{\pi(E_6), \pi(E'_6)\}) \, .$$
Moreover,~the bitangents
$\pi(E_1), \pi(E'_1), \ldots, \pi(E_6), \pi(E'_6)$
are distinct, as the intersection
number~$2$
is excluded, according to our~assumptions. Thus,~we indeed have a Steiner hexad.
}
\eop
\end{lem}

\begin{rems}
\label{revert}
\begin{iii}
\item
Assume~that
$S$
has a conic bundle with six split~fibres
$F_i = E_i \cup E'_i$,
for~$i = 1, \ldots, 6$.
Then~the assumptions of the criterion above are satisfied and hence
$\{(\pi(E_1),\pi(E'_1)), \ldots, (\pi(E_6),\pi(E'_6))\}$
is a Steiner~hexad.
\item
In~particular, we may revert the conclusion of Proposition~\ref{detf}.d). Indeed,~suppose that a double cover
of~$\Pb^2_k$,
ramified at the nonsingular
quartic~$C$,
has a
\mbox{$k$-rational}
conic bundle, with geometrically six degenerate fibres each splitting into two~lines. Then~every irreducible
component~$E$
is automatically an exceptional curve, as one has
$0 = EF = E(E+E') = E^2 + 1$,
for
$F$
the class of the~fibre. Therefore, the criterion applies and yields a Galois invariant Steiner hexad
for~$C$.
As~a consequence, there is a
\mbox{$k$-rational}
class
in~$\Pic(C)_2$,
as~well.
\end{iii}
\end{rems}

\begin{coro}
\label{int_Shex}
Let\/~$C \subset \Pb^2$
be a nonsingular plane quartic over an algebraically closed field
and\/~$L$
be a~bitangent. Then~the intersection of all Steiner hexads
containing\/~$L$
consists only
of\/~$L$.\medskip

\noindent
{\bf Proof.}
{\em
We~consider a double cover
of~$\Pb^2$,
ramified
at~$C$,
and work in the blown-up~model. Without~loss of generality, let us assume that
$L = \pi(e_1)$.
Then~Lemma \ref{crit_hex} yields the six Steiner hexads
$\Hb_2, \ldots, \Hb_7$
containing~$L$,
which are given~by
$$\Hb_k := \{\pi(e_i), \pi(l_{ik}) \mid 1 \le i \le 7, i \neq k\} \, ,$$
for~$k = 2, \ldots, 7$.
Note~here that the assumption
$i \neq k$
is necessary to have
$l_{ik}$
even defined. Furthermore,~it ensures that
$e_i \cdot l_{i'k} = 0$
for~$i \neq i'$.
Clearly,~$L = \pi(e_1)$
is the only bitangent
$\Hb_2, \ldots, \Hb_7$
all have in~common.%
}%
\eop
\end{coro}

\begin{prop}
\label{sympl}
Let\/~$C \subset \Pb^2$
be a nonsingular plane quartic over an algebraically closed field and\/
$S$
be the double cover
of\/~$\Pb^2$,
ramified
at\/~$C$.
We~equip\/
$\Pic(S)[\pi_*]/2\Pic(S)[\pi_*]$
with the\/
$\bbF_{\!2}$-valued
symplectic pairing, induced by the intersection pairing, and\/
$\Pic(C)_2$
with the Weil pairing.

\begin{abc}
\item
Then~the restriction
$$\overline{r}\colon \Pic(S)[\pi_*]/2\Pic(S)[\pi_*] \longrightarrow \Pic(C)_2$$
is a symplectic epimorphism.
\item
The kernel
of\/~$\overline{r}$
is the radical of\/
$\Pic(S)[\pi_*]/2\Pic(S)[\pi_*]$,
which is generated by the class of\/
$\smash{\calO_S(\widetilde{E}-E)}$,
for an arbitrary exceptional
curve\/~$E$
and\/~$\smash{\widetilde{E} := g(E)}$.
\end{abc}\smallskip

\noindent
{\bf Proof.}
{\em
First of all, let
$E_1, \ldots, E_7$
be mutually skew exceptional curves~and put
$\smash{\widetilde{E}_1 := g(E_1)}$.
As~seen in~\ref{Picardgr}.ii), up to the odd index three,
$\Pic(S)[\pi_*]$
is generated by the sheaves
$\calO_S(E_2-E_1)$,
\ldots,
$\calO_S(E_7-E_1)$,
and
$\smash{\calO_S(\widetilde{E}_1-E_1)}$.
Their~intersection matrix is,~obviously
$$
\left(
\begin{array}{rrrrrrr}
-2 & -1 & -1 & -1 & -1 & -1 & -2 \\
-1 & -2 & -1 & -1 & -1 & -1 & -2 \\
-1 & -1 & -2 & -1 & -1 & -1 & -2 \\
-1 & -1 & -1 & -2 & -1 & -1 & -2 \\
-1 & -1 & -1 & -1 & -2 & -1 & -2 \\
-1 & -1 & -1 & -1 & -1 & -2 & -2 \\
-2 & -2 & -2 & -2 & -2 & -2 & -6
\end{array}
\right) \, .
$$
Moreover,~the class of
$\smash{\calO_S(\widetilde{E}_1-E_1)}$
is clearly an element in the kernel
of~$\overline{r}$.

Thus,~one only has to show that, for
$v_i := r(\calO_S(E_i-E_1))$,
the Weil pairing
$\langle v_i, v_j\rangle$
evaluates nontrivially whenever
$i \neq j$,
for
$i,j \in \{2, \ldots, 7\}$.
In~order to do this, let us recall that the Riemann-Mumford relations \cite[Theorem 5.1.1]{Do} ensure that
$$\langle\varepsilon, \delta\rangle = h^0(C, \varepsilon+\delta+\eta) + h^0(C, \varepsilon+\eta) + h^0(C, \delta+\eta) + h^0(C, \eta) \pmod 2\, ,$$
for
$\eta \in \Pic(C)$
any invertible sheaf such that
$\smash{\eta^{\otimes 2} \cong \calO_{\Pb^2}(1)|_C}$.
Therefore,~taking
$\eta := r(\calO_S(E_1))$,
we~find
\begin{eqnarray*}
\langle v_i, v_j \rangle & = & h^0(C, r(\calO_S(E_i + E_j - E_1))) \\
 && \hspace{0.3cm} {} + h^0(C, r(\calO_S(E_i))) + h^0(C, r(\calO_S(E_j))) + h^0(C, r(\calO_S(E_1))) \pmod 2\, .
\end{eqnarray*}
But
$h^0(C, r(\calO_S(E_i))) = h^0(C, r(\calO_S(E_j))) = h^0(C, r(\calO_S(E_1))) = 1$,
as shown in (\ref{theta_odd}), so we only need to verify~that
$$h^0(C, r(\calO_S(E_i + E_j - E_1))) = 0 \, .$$
This is exactly Lemma \ref{lem_theta_even}~below.
\eop
}
\end{prop}

\begin{rems}
\begin{iii}
\item
The fact that
$\overline{r}$
is symplectic has been stated for the first time in \cite[Ch.~IX, Sec.~1, Lemma~2]{DO}. Our~proof follows the ideas of \cite[Remark~2.11]{Za} and is very different from the original~one.
\item
In~much more generality, A.\,N.\ Skorobogatov~\cite[Corollary~3.2.iii)]{Sk17} constructs a canonical homomorphism
$\Phi$
the other way round and shows that there is a short exact~sequence
$$0 \longrightarrow H^1_\et(C, \bbZ/2\bbZ) \stackrel{\Phi}{\longrightarrow} H^2_\et(S, \bbZ/2\bbZ) / \pi^* H^2_\et(\Pb^2, \bbZ/2\bbZ) \stackrel{\pi_*}{\longrightarrow} H^2_\et(\Pb^2, \bbZ/2\bbZ) \longrightarrow 0$$
of Galois~modules. The homomorphism
$\Phi$
is given by any section of the homomorphism
$H^2_\et(S, \bbZ/2\bbZ)[\pi_*] \to H^1_\et(C, \bbZ/2\bbZ)$
that is induced by the
restriction~$\overline{r}$
\cite[Lemma~3.3]{Sk17}.
It~is, unfortunately, not discussed in~\cite{Sk17} whether or in which generality
$\Phi$
is~symplectic.
\end{iii}
\end{rems}

\begin{lem}
\label{lem_theta_even}
Let\/~$C \subset \Pb^2$
be a nonsingular plane quartic over an algebraically closed field and\/
$S$
the double cover
of\/~$\Pb^2$,
ramified
at\/~$C$.
Then,~for three mutually skew exceptional curves\/
$E$,
$E'$,
and\/~$E''$
on\/~$S$,
one has
\begin{equation}
\label{theta_even}
h^0(C, r(\calO_S(E + E' - E''))) = 0 \, .
\end{equation}
{\bf Proof.}
{\em
The~only way for
$r(\calO_S(E + E' - E''))$
to have a section is that there is a fourth bitangent
$L$
of~$C$
such that
$\Pi(\{L, \pi(E)\}) = \Pi(\{\pi(E'), \pi(E'')\})$.
Indeed,~in a notation analogous to that used before Remark~\ref{bita_div}, one has
$$\calL := r(\calO_S(E' - E'')) = \calO_C((p'_1)+(p'_2)-(p''_1)-(p''_2)) \in \Pic(C)_2 \, ,$$
and a section of
$r(\calO_S(E + E' - E''))$
would yield two points
$p'''_1$
and
$p'''_2 \in C$
such that
$\calO_C((p'''_1)+(p'''_2)-(p_1)-(p_2)) = \calL$,
too, for
$p_1$
and~$p_2$
the points of tangency of
$\pi(E)$
with~$C$.
In~particular, the~divisor
$$D := 2(p'''_1)+2(p'''_2)-2(p_1)-2(p_2) \in \Div(C)$$
would be~principal.
However,~$\smash{\Gamma(C, \calO_C(2(p_1) + 2(p_2))) = \langle \frac{T_0}s, \frac{T_1}s, \frac{T_2}s \rangle}$,
for
$s$
the linear form defining the bitangent
$\pi(E)$,
according to formula~(\ref{lin_schnitt}). In~particular, one must have
$D = \div(l/s)$
for a further linear
form~$l$.
But~then,
$L := Z(l)$
is clearly a bitangent touching
$C$
in
$p'''_1$
and~$p'''_2$.

Take an exceptional curve
$E'''$
on~$S$
such that
$\pi(E''') = L$.
As~$r'$
induces an epimorphism
$\overline{r}\colon \Pic(S)[\pi_*]/2\Pic(S)[\pi_*] \longrightarrow \Pic(C)_2$
with kernel generated by the class
of~$\smash{\calO_S(\widetilde{E}-E)}$,
we find that
\begin{eqnarray*}
\calO_S(E''' - E) & \equiv & \calO_S(E' - E'')\hspace{0.83cm} \pmod {\langle2\Pic(S)[\pi_*], \calO_S(\widetilde{E}-E) \rangle} \, , \qquad\text{i.e.\,,} \\
\calO_S(E''') & \equiv & \calO_S(E+E'+E'') \pmod {\langle2\Pic(S)[\pi_*], \calO_S(\widetilde{E}-E) \rangle} \, .
\end{eqnarray*}

In~order to complete the argument, let us work in the blown-up model. One~may assume without restriction that
$E=e_1$,
$E'=e_2$,
and
$E''=e_3$.
Then
$E'''$
must have odd intersection number
($-1$
or
$1$)
with
$e_1$,
$e_2$,
and
$e_3$
and even intersection number
($0$~or~$2$)
with
$e_4,\ldots,e_7$,
or the other way~round. A~short consideration shows, however, that such an exceptional curve does not~exist.
}
\eop
\end{lem}

\begin{rem}
Formula (\ref{theta_even}) is classically known, as
$\pi(E) = \pi(e_1)$,
$\pi(E') = \pi(e_2)$,
and
$\pi(E'') = \pi(e_3)$
form part of the so-called Aronhold set
$\{\pi(e_1), \ldots, \pi(e_7)\}$,
cf.~\cite[Section 6.1.2]{Do}.
\end{rem}

\begin{coro}[The admissible subgroup in degree~28 versus that in degree~56]
Let\/~$C$
be a nonsingular plane quartic over an algebraically closed field and\/
$S$
be the double cover
of\/~$\Pb^2$,
ramified
at\/~$C$.
Write\/~$\smash{\widetilde\Omega_C}$
for the set of the 56 exceptional curves and\/
$\smash{\widetilde{G} \subset \Sym \widetilde\Omega_C \cong S_{56}}$
for the maximal subgroup that respects the intersection pairing
on\/~$S$.
Moreover,~put
\begin{align*}
G := \{\sigma \in \Sym \Omega_C \mid{} & \text{The permutation\/ } \sigma \text{ of the 28 bitangents respects the} \\[-1.7mm]
 & \text{\,Steiner hexads, and the operation on the 63 Steiner hexads} \\[-1.7mm]
 & \hspace{1.5cm} \text{induces a symplectic automorphism of\/ } \smash{H^1_\et(C, \bbZ/2\bbZ)} \} \, .
\end{align*}%
Then

\begin{abc}
\item
$\smash{\widetilde{G} \cong W(E_7)}$
and\/
$G \cong \Sp_6(\bbF_{\!2})$.
\item
The group\/
$\smash{\widetilde{G}}$
has the pairs of exceptional curves, lying over the same bitangent, as a block~system. The permutation representation\/
$\iota\colon W(E_7)/Z \to \Sym \Omega_C \cong S_{28}$
induced on the blocks is faithful.
\item
The image
of\/~$\iota$
is
exactly\/~$G$.
\end{abc}\smallskip

\noindent
{\bf Proof.}
{\em
First,~the canonical homomorphism
$\smash{G \to \Aut(H^1_\et(C, \bbZ/2\bbZ)) \cong \Sp_6(\bbF_{\!2})}$
is~injective. Indeed,~for a permutation
$\sigma \in G \subset \Sym \Omega_C \cong S_{28}$
to lie in the kernel,~the induced operation
$\smash{\sigma^{(2)} \in \Sym \Omega_C^{(2)} \cong S_{\frac{28\cdot27}2}}$
on 2-sets must keep all Steiner hexads in~place. According~to Corollary~\ref{int_Shex}, this yields
$\smash{\sigma^{(2)} = \id}$,
which suffices
for~$\sigma = \id$.

On~the other hand, the fact that
$\smash{\widetilde{G} \cong W(E_7)}$
is shown in \cite[Theorem~23.9]{Ma}. The~operation
of~$\smash{\widetilde{G}}$
respects the pairs
$\smash{\{E, \widetilde{E}\}}$
as these are the only ones with intersection
number~$2$.
Moreover,~the centre, which is generated by the Geiser involution, clearly fixes all blocks, such that there is an induced permutation representation
$\iota\colon W(E_7)/Z \to \Sym \Omega_C$.

A~nontrivial element in kernel
of~$\iota$
would be represented by some
$\tau \in W(E_7)$
that flips some but not all of the
pairs~$\smash{\{E, \widetilde{E}\}}$.
The~projection formula shows, however, that, for arbitrary exceptional curves
$E$
and~$E'$,
one has
\begin{equation}
\label{projf}
\smash{EE' + E\widetilde{E}' = E(E' + \widetilde{E}') = E \cdot \pi^*L = \pi_*E \cdot L = 1} \, ,
\end{equation}
for
$\smash{\widetilde{E}' = g(E')}$.
The~assumption that
$\tau$
would act in the way that
$E \mapsto E$,
$\smash{E' \mapsto \widetilde{E}'}$,
and
$\smash{\widetilde{E}' \mapsto E'}$
is hence contradictory, as one has, without restriction,
$\smash{EE' = 0}$
and
$\smash{E\widetilde{E}' = 1}$.
Thus,~we see that
$\iota$
is injective, which completes the proof of~b).

Finally,~the operation of the group
$W(E_7)$
on~$\Pic(S)$,
and hence on
$\smash{\Pic(S)[\pi_*]}$,
respects the intersection pairing. Consequently,~the
$\bbF_{\!2}$-valued
symplectic pairing on
$\Pic(S)[\pi_*]/2\Pic(S)[\pi_*]$
is respected,~too. Therefore,~Proposition~\ref{sympl}.a) yields that 
$W(E_7)$
operates on
$\smash{\Pic(C)_2 \cong H^1_\et(C, \bbZ/2\bbZ)}$
via symplectic~automorphisms. In~particular, the homomorphism
$\iota\colon W(E_7)/Z \to \Sym \Omega_C$
factors
via~$G$,
and the following diagram~commutes,
$$
\xymatrix{
\Sym \Omega_C \cong S_{28} \ar@{=}[rr] && \Sym \Omega_C \cong S_{28} & \\
W(E_7)/Z \ar@{^{(}->}[u] \ar@{->}[rr]^{\;\;\;\;\;\;\;\;\;\overline{\iota}} && G \ar@{^{(}->}[r] \ar@{^{(}->}[u] & \Sp_6(\bbF_{\!2}) \, .
}
$$
Altogether, we have two injections
$W(E_7)/Z \hookrightarrow G \hookrightarrow \Sp_6(\bbF_{\!2})$
in a row. As~both groups,
$W(E_7)/Z$
and~$\Sp_6(\bbF_{\!2})$,
are of order
$1\,451\,520$,
the injections must be~bijective. In~other words, c)~is shown and the proof of~a) is~complete.
}
\eop
\end{coro}

\begin{coro}
Define the subgroup\/
$U_{63} \subset G$
as the stabiliser of a Steiner~hexad. Then

\begin{abc}
\item
The subgroup\/
$U_{63}$
is uniquely determined up to conjugation
in\/~$G$
and of
index\/~$63$.
As~a permutation group in degree 28,
$U_{63}$
is intransitive of orbit
type\/~$[12,16]$.
\item
As~an abstract group,
$U_{63}$
is isomorphic
to\/~$(S_2 \wr S_6) \cap A_{12}$,
operating on the size twelve orbit in the obvious~way. In~particular, the operation of an element\/
$\sigma \in U_{63}$
on the size twelve orbit determines\/
$\sigma$~completely.
\end{abc}\smallskip

\noindent
{\bf Proof.}
{\em
Let~us consider a non-singular model quartic
$C \subset \Pb^2$
over an algebraically closed field and put
$S$
to be the double cover
of~$\Pb^2$,
ramified
at~$C$.\smallskip

\noindent
a)
According~to Definition~\ref{St_hex}, Steiner hexads are in bijection with the 63 nonzero elements in
$\Pic(C)_2$.
Moreover,~the
group~$G \cong \Sp_6(\bbF_{\!2})$
operates transitively on those nonzero~elements and therefore as well on the 63 Steiner~hexads. This implies the two first assertions. 

The~final one is easily checked by an experiment in {\tt magma}~\cite{BCP}, which shows the~following. The,~up to conjugation, unique
\mbox{index-$63$}
subgroup of the only simple permutation group of order
$1\,451\,520$
in
degree~$28$
has orbit type
$[12, 16]$.\smallskip

\noindent
b)
Once again, the subgroup
$U_{63} \subset G$
distinguishes a Steiner hexad and hence a nonzero element 
of~$\Pic(C)_2$.
According~to Proposition~\ref{detf}.d),
$S$
has an associated conic~bundle. As~before, let us write
$E_1, E'_1, \ldots, E_6, E'_6$
for the irreducible components of its six split~fibres. Then~the distinguished Steiner hexad is simply
$\{(\pi(E_1), \pi(E'_1)), \ldots, (\pi(E_6), \pi(E'_6))\}$.
Having~numbered these twelve bitangents from
$1$
to~$12$,
the operation
of~$U_{63}$
defines a homomorphism
$i\colon U_{63} \to S_{12}$.\smallskip

\noindent
{\em Injectivity.}
Assume that
$\sigma \in U_{63}$
fixes each of the twelve~bitangents. Then~a lift
$\smash{\widetilde\sigma \in \widetilde{G}}$
of~$\sigma$
either fixes
$E_1, E'_1, \ldots, E_6, E'_6$
or sends each of these exceptional curves to its image under the Geiser involution. Without~restriction, we may assume that
$\smash{\widetilde\sigma}$
fixes the~curves.

Next,~we observe that the sum over all 56 exceptional curves
in~$\Pic(S)$
is~$-28[K]$.
This~shows that any permutation
in~$\smash{\Sym \widetilde\Omega_C \cong S_{56}}$
fixes the canonical
class~$[K]$.
Moreover,~as seen in~(\ref{rank8}), together
with~$[K]$,
the curves
$E_1, E'_1, \ldots, E_6, E'_6$
generate
$\Pic(S)$
up to a finite~index.
Consequently,~$\smash{\widetilde\sigma}$
operates identically on the whole
of~$\Pic(S)$.
It~must therefore fix each of the 56 exceptional curves, which shows
that~$\smash{\widetilde\sigma = \id}$.\smallskip

\noindent
{\em Image.}
The~maximal permutation group that respects the intersection matrix of
$E_1, E'_1, \ldots, E_6, E'_6$
is clearly the wreath product
$S_2 \wr S_6$.
However,~the Faddeev reciprocity law for conic bundles~\cite[Corollary~3.3]{Sk13} implies that Galois operations permuting these twelve curves are limited to even permutations. We~claim that, for a model surface, the action of
$U_{63}$
is provided by a Galois operation. Indeed,~the results of T.~Shioda~\cite[Section~3]{Shi} and R.~Ern\'e~\cite[Corollary~3]{Er} provide a nonsingular plane quartic
over~$\bbQ$
of the kind that
$\smash{\Gal(\overline\bbQ/\bbQ)}$
operates on the bitangents
via~$G$.
Hence,~for a certain number field
$l_{63}$
of
degree~$63$,
$\smash{\Gal(\overline\bbQ/l_{63})}$
acts precisely via
$U_{63} \subset G$.
Therefore,~the image of
$i\colon U_{63} \to S_{12}$
is contained in
$(S_2 \wr S_6) \cap A_{12}$.
As~both groups,
$U_{63}$
and~$(S_2 \wr S_6) \cap A_{12}$,
are of order
$23\,040$,
equality~holds.
}
\eop
\end{coro}

\section{Conic bundles with prescribed Galois operation I--\\The construction}
\label{drei}

Let~a
field~$k$,
a Galois extension
$l/k$,
and an injection
$g := \Gal(l/k) \stackrel{i}{\hookrightarrow} S_2 \wr S_6$
be~given. Then~one has the~subgroups
$$B := \im(\pro \!\circ i\colon g \to S_6) \subseteq S_6
\quad{\rm and}\quad
H := \ker(\pro \!\circ i\colon g \to S_6) \subseteq g \, .$$
The~normal subgroup
$H \subseteq g$
defines an intermediate
field~$k'$
such that
$\Gal(k'/k) \cong B$
and
$\Gal(l/k') = H$.
As~$H \subseteq (\bbZ/2\bbZ)^6$,
the latter is a Kummer extension.
I.e.,~one has
$\smash{l = k'(\sqrt{A})}$,
for a subgroup
$A \subseteq k'^*/(k'^*)^2$.
The~lemmata below analyse this situation from the field-theoretic point of~view.

For~this, let us note that the group
$\smash{g \stackrel{i}{\hookrightarrow} S_2 \wr S_6}$
naturally operates on twelve objects
$1_a, 1_b, \ldots, 6_a, 6_b$
forming six~pairs. For~every
$n \in \{1, \ldots, 6\}$,
the stabiliser
$\Stab_g(n_a) = \Stab_g(n_b)$
is a subgroup of index
$1$
or~$2$
in~$\Stab_g(\{n_a, n_b\})$,
hence~normal.

Let~us denote the subfield
of~$l$,
corresponding
to~$\Stab_g(\{n_a, n_b\})$
under the Galois correspondence,
by~$k_n$.
Similarly,~we write
$l_n$
for the subfield corresponding
to~$\Stab_g(n_a)$.
According~to these definitions, one has
$k_n \subseteq k'$
and
$k_n \subseteq l_n$,
the latter extension being of degree
$1$
or~$2$.
Moreover,~$k' = k_1 \cdots k_6$
and
$l = l_1 \cdots l_6$.

\begin{sub}
\label{tech}
Let\/~$k$
be a field,
$d > 0$
be an integer, and\/
$l$
an \'etale
\mbox{$k$-algebra\/}
of
degree\/~$d$.

\begin{abc}
\item
Then,~for every\/
$x \in l^*$,
there is a dominant morphism\/
$q_x\colon \Ab(l) \to \Ab(l)$
of\/
\mbox{$k$-schemes}
that induces on\/
\mbox{$k$-rational}
points the mapping\/
$l \to l$,
$t \mapsto xt^2$.
\item
There is a dominant morphism\/
$c_l\colon \Ab(l) \to \Ab(k[T]_d)$
of\/~\mbox{$k$-schemes}
inducing on\/
\mbox{$k$-rational}
points the mapping\/
$l \to k[T]_d$,
$t \mapsto \chi_t$,
for\/
$\chi_t$
the characteristic polynomial of the multiplication map\/
$\cdot t\colon l \to l$.
\item
Let\/
$d_1$
and\/~$d_2$
be positive integers such that\/
$d_1 + d_2 = d$.
Then~there is a dominant morphism\/
$m_{d_1,d_2}\colon \Ab(k[T]_{d_1}) \times \Ab(k[T]_{d_2}) \to \Ab(k[T]_d)$
of\/~\mbox{$k$-schemes}
that induces on\/
\mbox{$k$-rational}
points the mapping\/
$k[T]_{d_1} \times k[T]_{d_2} \to k[T]_d$,
$(f_1, f_2) \mapsto f_1f_2$.
\end{abc}\smallskip

\noindent
{\bf Proof.}
{\em The mappings given are clearly polynomial, so they define morphisms of
\mbox{$k$-schemes}.

Moreover,~dominance may be tested after base extension to the algebraic~closure.
Then,~$l \cong k^d$.
Thus,~by inspecting
\mbox{$k$-rational}
points, one finds that all three types of morphisms are always quasi-finite, with generic fibres of sizes
$2^d$,
$d!$,
and~$\smash{\binom{d }{d_1}}$,
respectively. As~source and target are of the same dimension, this suffices for~dominance.
}
\eop
\end{sub}

\begin{lem}
Let~a Galois extension\/
$l/k$
and an injection\/
$\Gal(l/k) \hookrightarrow S_2 \wr S_6$
be~given, and let\/
$A$,
$B$,
and\/~$k'$
have the same meaning as~above.

\begin{abc}
\item
Then~there is a natural ordered set\/
$(\alpha_1, \ldots, \alpha_6)$
contained in\/
$A \subseteq k'^*/(k'^*)^2$
such~that
\begin{iii}
\item
The~elements\/
$\alpha_1, \ldots, \alpha_6$
generate the abelian
group\/~$A$.
\item
The~Galois group\/
$\Gal(k'/k) \cong B$
permutes\/
$\alpha_1, \ldots, \alpha_6 \in k'^*/(k'^*)^2$
exactly as described by the natural inclusion\/
$B \subseteq S_6$.
\end{iii}
\item
There exist lifts\/
$A_1, \ldots, A_6 \in k'^*$
of\/
$\alpha_1, \ldots, \alpha_6$
that are still permuted by\/
$\Gal(k'/k)$.
\end{abc}\smallskip

\noindent
{\bf Proof.}
{\em
Let~us note that
$k'^*/(k'^*)^2$
is naturally acted upon
by~$\Gal(k'/k)$.
Indeed,~an automorphism
of~$k'$
sends
$(k'^*)^2$
to~itself.\smallskip

\noindent
a)
The~group
$\Gal(l/k') = \Gal(k'(\sqrt{A})/k') = H$
is abelian of
exponent~$2$.
Hence,~according to Kummer theory,
$H$
is the dual of the
group~$A \subseteq k'^*/(k'^*)^2$
and, consequently,
$A \cong H^\vee$
holds as~well. Furthermore,~the standard linear forms on
$H \subseteq (\bbZ/2\bbZ)^6$,
$$e_i: (z_1, \ldots, z_6) \mapsto z_i \, ,$$
for
$i = 1, \ldots, 6$,
generate
$H^\vee$
and are permuted by
$\Gal(k'/k) \cong B \subseteq S_6$
in the natural~manner. Thus,~the elements
$\alpha_1, \ldots, \alpha_6 \in A \cong H^\vee$,
corresponding to these linear forms, satisfy conditions~i) and~ii).\smallskip

\noindent
b)
The~assertion simply means that each
$\alpha_n$
may be lifted
to~$k'^*$
in such a way that the size of its
\mbox{$g$-orbit}
remains the~same. I.e.,~such that
$\Stab_g(A_n) = \Stab_g(\alpha_n)$.
Indeed,~if
$B \subseteq S_6$
is transitive then one may simply lift
$\alpha_1$
in this way, and take the
\mbox{$B=\Gal(k'/k)$-orbit}
of the~lift. In~the intransitive case, the same has to be done separately for each orbit
in~$\{1, \ldots, 6\}$.

In~order to verify liftability of this particular kind for one of
the~$\alpha_n$,
recall that the classes
$\alpha_1, \ldots, \alpha_6 \in k'^*/(k'^*)^2$
are characterised by the property that
$\smash{\pm\sqrt{\mathstrut\alpha_1}}$,
\ldots, $\smash{\pm\sqrt{\mathstrut\alpha_6}}$
are acted upon
by~$g$
in the same way as the symbols
$1_a, 1_b, \ldots, 6_a, 6_b$.
Thus,~we need elements
$\sqrt{A_n} \in l$
of the kind~that
\begin{equation}
\label{orbit_size}
\Stab_g(\sqrt{A_n}) = \Stab_g(n_a) \,.
\end{equation}
Then~clearly
$\smash{\Stab_g(A_n) = \Stab_g(\{n_a, n_b\}) = \Stab_g(\alpha_n)}$.

For~this, we know that
$l_n$~is
an at most quadratic extension
of~$k_n$.
Hence,
$\smash{l_n = k_n(\sqrt{A_n})}$,
for some
$A_n \in k_n^*$.
As~the fields
$k_n$
are permuted
by~$\Gal(k'/k)$
in the same way as the sets
$\{n_a, n_b\}$,
at least as long as it is a primitive element of the
field~$k_n$,
$A_n$~has
the property~(\ref{orbit_size})~required.
}
\eop
\end{lem}

\begin{lem}
\label{pol_12}
Let~a Galois extension\/
$l/k$
and an injection\/
$\Gal(l/k) \hookrightarrow S_2 \wr S_6$
be given, and let\/
$\alpha_n$,
$A_n$,
and\/~$k'$
be as~above. Choose~lifts\/
$A_1, \ldots, A_6 \in k'^*$
of\/
$\alpha_1, \ldots, \alpha_6 \in k'^*/(k'^*)^2$
that form a\/
$\Gal(k'/k)$-invariant
set.

\begin{abc}
\item
Then the polynomial\/
$F \in k[T]$
such~that
$$F(T) := (T-A_1) \cdots (T-A_6)$$
has splitting
field\/~$k'$.
\item
Furthermore,~$k(\sqrt{A_1}, ..., \sqrt{A_6}) = l$.
I.e.,~$l$
is the splitting field
of\/~$F(T^2)$.
Finally,~the operation of\/
$\Gal(l/k)$
on the roots\/
$\pm\sqrt{A_i}$
agrees with the natural operation of\/
$S_2 \wr S_6$
on twelve objects forming six~pairs.
\end{abc}\smallskip

\noindent
{\bf Proof.}
{\em
a)
By~construction, the polynomial
$F \in k'[T]$
is
$\Gal(k'/k)$-invariant,
hence
$F \in k[T]$.
On~the other hand, the splitting field
of~$F$
is
$k(A_1, \ldots, A_6)$,
which is clearly contained
in~$k'$.
The~claim follows, since the Galois group
of~$F$
is
$B \subseteq A_6$
and hence coincides with that
of~$k'$.\smallskip

\noindent
b)
The splitting field of
$F(T^2)$
is
$k(\sqrt{A_1}, ..., \sqrt{A_6})$,
which is equal
to~$l$,
according to our~construction. The~final assertion is~obvious.
}
\eop
\end{lem}

\begin{coro}
\label{Zariski}
Let~a Galois extension\/
$l/k$
and an injection\/
$\Gal(l/k) \hookrightarrow S_2 \wr S_6$
be given, where the
field\/~$k$
is~infinite. Then~the set of all polynomials\/
$F$
as in Lemma \ref{pol_12} is Zariski dense
in\/~$\Ab(k[T]_6)$.\medskip

\noindent
{\bf Proof.}
{\em
Let
$O_1, \ldots, O_m \subseteq \{1,\ldots,6\}$
be the orbits under the operation
of~$B \subseteq S_6$.
For~each of them, let us choose a representative
$n_i \in O_i$
and a lift
$\smash{A_{n_i}^{(0)} \in k_{n_i}^*}$.
Then~the polynomials in the sense above certainly include those of~type
$$\prod_{i=1}^m \chi_{A_{n_i}^{(0)} \cdot t_i^2} (T) \, ,$$
for all
$\smash{t_i \in k_{n_i}^*}$
such that
$A_{n_i}^{(0)} \!\cdot\! t_i^2 \in k_{n_i}$
is a primitive~element.

In~order to prove Zariski density, let,~at first,
$i \in \{1, \ldots, m\}$
be~arbitrary. Then,
as~$k_{n_i} \supseteq k$
is an infinite field,
$\smash{k_{n_i} = \Ab(k_{n_i})(k)}$
is Zariski dense
in~$\Ab(k_{n_i})$.
Thus, Sublemma~\ref{tech}.a) shows that the elements
$$q_{A_{n_i}^{(0)}}(t_i) = A_{n_i}^{(0)} \!\cdot\! t_i^2 \in k_{n_i} = \Ab(k_{n_i})(k) \, ,$$
for~$t_i \in k_{n_i}^*$,
form a Zariski dense subset
in~$\Ab(k_{n_i})$,~too.
The~same is still true for the subset of these elements that are~primitive. Indeed,~being a primitive element is a Zariski open~condition. Consequently,~by Sublemma \ref{tech}.b), the polynomials
$$c_{k_{n_i}}(A_{n_i}^{(0)} \!\cdot\! t_i^2) = \chi_{A_{n_i}^{(0)} \cdot t_i^2} (T) \in k[T]_{\#O_{n_i}} = \Ab(k[T]_{\#O_{n_i}})(k)$$
are Zariski dense
in~$\Ab(k[T]_{\#O_{n_i}})$.
Part~c) of Sublemma~\ref{tech} finally implies the~claim.
}
\eop
\end{coro}

\begin{prop}[A conic bundle with six split fibres and prescribed Galois action]
\label{cb_gal}
Let\/~$k$
be a field and\/
$F(T) = T^6 + a_5T^5 + \cdots + a_1T + a_0 \in k[T]$
be a separable, monic polynomial of degree six, where\/
$a_0 = c^2$,
for
some\/~$c \in k^*$.

\begin{abc}
\item
Then~there exist exactly two pairs\/
$(g,h)$
of binary quadratic forms such~that
$$\det M(1,T) = -F(T) \,,$$
for
\begin{equation}
\label{matrix}
M(s,t) :=
\left(
\begin{array}{ccc}
-st + t^2 &  st & g(s,t) \\
 st       & s^2 & t^2    \\
 g(s,t)   & t^2 & h(s,t)
\end{array}
\right) .
\end{equation}
For~these, one
has~$g(1,0) = \pm c$.
\item
Let\/~$(g,h)$
be one of the pairs as in~a), that
with\/~$g(1,0) = c$.
Then~the equation
\begin{equation}
\label{cb}
(T_0\, T_1\, T_2)\, M(s,t)
\!\left(
\begin{array}{c}
T_0 \\ T_1 \\ T_2
\end{array}
\right)
= 0
\end{equation}
defines a
hyperplane\/~$B_{F,c}$
of
bidegree\/~$(2,2)$
in\/~$\Pb^1 \times \Pb^2$.
\begin{iii}
\item
The generic fibre of the projection
$\pro_1\colon B_{F,c} \to \Pb^1$
is a~conic. Degenerate~fibres occur exactly over the roots
of\/~$F$.
\item
Each of the six reducible fibres geometrically splits into two~lines. The irreducible components are defined over the splitting
field\/~$l$
of\/~$F(T^2)$
and permuted
by\/~$\Gal(l/k)$
in the same way as the roots
of\/~$F(T^2)$.
\item
The projection\/
$\pro_2\colon B_{F,c} \to \Pb^2$
is a double cover
of\/~$\Pb^2$,
ramified at a quartic
curve~$C_{F,c}$.
For~a generic choice of the
polynomial\/~$F$,
the curve
$C_{F,c}$~is
nonsingular.
\end{iii}
\end{abc}\smallskip

\noindent
{\bf Proof.}
{\em
a)
A direct calculation shows that the determinant
of~$M(1,T)$
is~exactly
\begin{equation}
\label{detconst}
-Th(1,T) + 2T^3g(1,T) - g^2(1,T) - (T^2-T)T^4 = -Th(1,T) + T^5 - [T^3 - g(1,T)]^2 \, .
\end{equation}
Writing~$g(1,T) = g_2 T^2 + g_1 T + g_0$,
this term takes the~form
$$-Th(1,T) - T^6 + (2g_2 +1)T^5 - (g_2^2-2g_1)T^4 + \cdots - g_0^2 \, ,$$
so
$F$
determines the coefficients
$\smash{g_2 := -\frac{a_5+1}2}$
and~$\smash{g_1 := \frac{-a_4+g_2^2}2}$
uniquely and the absolute term
$\smash{g_0 := \pm\sqrt{\mathstrut a_0} = \pm\sqrt{\mathstrut c^2}}$
up to~sign. In~particular,
$g(1,0) = g_0 = \pm{c}$.
\mbox{Finally,~for} either choice
of~$g$,
the coefficients
of~(\ref{detconst})
at~$T^i$,
for~$i = 1$,
$2$,
and~$3$,
uniquely determine the quadratic polynomial
$h(1,T)$
to~be
$$h(1,T) := \frac{F - [T^3-g(1,T)]^2 + T^5}T \, .$$

\noindent
b.i)
Split fibres correspond to zeroes of the determinant
of~$M(s,t)$.\smallskip

\noindent
ii)
The~first claim is that
$\rk M(s,t) = 2$
for each
$(s\!:\!t) \in \Pb^1$
such that
$\det M(s,t) = 0$.
For~this, we observe that the upper left
$2\times2$-minor
of~$M(s,t)$
is
$(-s^3t)$,
which vanishes only at
$0$
and~$\infty$.
However,~all split fibres occur over points on the affine line and
$0$
is excluded, since
$F(0) = a_0 = c^2 \neq 0$.

Finally,~let
$t_0$
be any root
of~$F$.
Then~the fibre
$(B_{F,c})_{(1:t_0)}$
is the degenerate conic defined
over~$k(t_0)$
that is given by the
$3\times3$-matrix
$M(1,t_0)$.
As~$M(1,t_0)$
has vanishing determinant and a principal
$2\times2$-minor
of~$(-t_0)$,
the fibre
$\smash{(B_{F,c})_{(1:t_0)}}$
is projectively equivalent to the degenerate conic, given by
$T_0^2 - t_0 T_1^2 = 0$.
This~latter conic clearly splits into two lines
over~$\smash{k(\sqrt{\mathstrut t_0})}$.

For~$t_1,\ldots,t_6$
the roots
of~$F$,
the field of definition of the twelve irreducible components is hence
$\smash{k(\sqrt{\mathstrut t_1}, \ldots, \sqrt{\mathstrut t_6})}$,
which is exactly the splitting field
of~$F(T^2)$.
The~final claim is~obvious.\smallskip

\noindent
iii)
The left hand side of equation~(\ref{cb}) may be considered as a binary quadratic form, the coefficients of which are quadratic forms in
$T_0$,
$T_1$,
and~$T_2$,
i.e.~as
$$Q_0(T_0,T_1,T_2)s^2 + Q_1(T_0,T_1,T_2)st + Q_2(T_0,T_1,T_2)t^2 \, .$$
As~such, its discriminant is the ternary quartic form
$\bm{Q}_{F,c} := Q_1^2 - 4Q_0Q_2$.
Thus,~$C_{F,c}$
is the vanishing locus
of~$\bm{Q}_{F,c}$. 

Using~{\tt magma}~\cite{BCP}, it is easy to write down
$\bm{Q}_{F,c}$
explicitly, in terms of the parameters
$a_1$,
\ldots,
$a_5$,
and~$c$.
We~used the first author's code~\cite{El} for calculating the Dixmier--\discretionary{}{}{}Ohno invariants of ternary quartic forms and found that the discriminant
of~$\bm{Q}_{F,c}$
is not the zero~polynomial. 
}
\eop
\end{prop}

\begin{rems}
\label{discr}
\begin{iii}
\item
The quartic form
$\bm{Q}_{F,c}$
occurs to us in the symmetric determinantal
form~$Q_1^2 - 4Q_0Q_2$.
This~is related to the fact that 
$C_{F,c}$
has a distinguished Steiner hexad, respectively a distinguished element in
$\Pic(C_{F,c})_2 \setminus \{\calO_{C_{F,c}}\}$,
cf.~Proposition~\ref{detf}.a). The quadratic forms
$Q_0$~and~$Q_2$
define particular contact conics
of~$C_{F,c}$.
\item
One may adopt the relative point of view, according to which formula~(\ref{cb}) actually describes a family
$\kappa\colon \bm{C} \to \Spec \bbQ[c,a_1,\ldots,a_5] = \Ab^6_\bbQ$
of plane~quartics.

Then~the discriminant
of~$\bm{Q}_{F,c}$
is a polynomial
in~$\bbQ[c,a_1,\ldots,a_5]$
that splits into two factors. One~of them is exactly the discriminant of the degree six
polynomial~$F$.
This~coincidence is, of course, not at all surprising, since multiple zeroes
of~$F$
cause a degenerate conic~bundle.

The~other factor, entering the discriminant
of~$\bm{Q}_{F,c}$
quadratically, reflects the fact that the total
scheme~$\bm{C}$
of the family
$\kappa$
is~singular. The~image
$\Delta_\sing := \kappa(\bm{C}_\sing)$
of the singular locus under projection to the parameter scheme
$\Ab^6$
is the divisor defined by this~factor. The~corresponding fibres generically have only one singular point that is an ordinary double~point.
\end{iii}
\end{rems}

\begin{rem}
The discriminant of the polynomial
$F(T^2)$~is,
up to square factors,
$(-1)^{\deg F} F(0)$.
In~particular, in our situation, the discriminant is a square
in~$k$,
such that one automatically has an~injection
$$i\colon \Gal(l/k) \hookrightarrow (S_2 \wr S_6) \cap A_{12} \, .$$

Indeed,~the general formula for the discriminants in a tower of fields (cf.~\cite[Corollary~III.2.10]{Ne}) shows
$$\Delta_{F(T^2)} = \Delta_F^2 \cdot N_{(k[T]/(F))/k}(r) \, ,$$
for~$r$
a root
of~$F$.
Therefore,
$\smash{\Delta_{F(T^2)} = \Delta_F^2 \cdot (-1)^{\deg F} F(0)}$.
But,~in our case, one has
$(-1)^{\deg F} F(0) = F(0) = a_0 = c^2$,
which shows the~claim.
\end{rem}

\section{Conic bundles with prescribed Galois operation II--\\The monodromy}
\label{vier}

\begin{ttt}
Let~an injection
$i\colon \Gal(l/k) \hookrightarrow (S_2 \wr S_6) \cap A_{12}$
be given and
$F$~be
a polynomial of the kind that
$l$
is the splitting field
of~$F(T^2)$.
Then~Proposition~\ref{cb_gal}.b.ii) yields a numbering of the split fibres of the conic bundle
surface~$B_{F,c}$
from
$1$
to~$6$,
such that
$\Gal(l/k)$
operates on the twelve irreducible components as described
by~$i$.
This~shows that, when a numbering is given on the split fibres, then the operation
of~$\Gal(l/k)$
on the twelve irreducible components agrees with that described
by~$i$
only up to an inner automorphism
of~$S_2 \wr S_6$.

According~to Lemma~\ref{crit_hex} and Remark~\ref{revert}.i), the projections of the six reducible fibres are twelve bitangents
of~$C_{F,c}$
forming a distinguished Steiner
hexad~$\Hb$.
Moreover,~this Steiner hexad is clearly Galois~invariant. When~the 28~bitangents
of~$C_{F,c}$
are equipped with a numbering, distinguishing the Steiner
hexad~$\Hb$,
then the above considerations apply. They~show that the operation
of~$\Gal(l/k)$
on the twelve bitangents
forming~$\Hb$
agrees with the one described
by~$i$,
but only up to conjugation by an element
of~$S_2 \wr S_6$.
\end{ttt}

\begin{rems}
\begin{iii}
\item
This~is not entirely adequate to our situation, as only subgroups of
$(S_2 \wr S_6) \cap A_{12}$
may~occur. Even~worse, as discussed in Remark~\ref{mark}.i), two numberings on the bitangents of a plane quartic that distinguish the same Steiner hexad may differ only by the conjugation with an element of
$U_{63} \cong (S_2 \wr S_6) \cap A_{12}$.
\item
Thus,~it may happen that Proposition~\ref{cb_gal} provides a quartic that behaves in an unwanted~way. The~operation
of~$\Gal(l/k)$
on the twelve bitangents
of~$C_{F,c}$,
forming the distinguished Steiner
hexad~$\Hb$,
might differ from the desired one, described
by~$i$.
The~difference would then be given by the outer
automorphism~$\Psi$
of~$(S_2 \wr S_6) \cap A_{12}$,
provided by the conjugation with an element from
$(S_2 \wr S_6) \setminus A_{12}$.
Fortunately,~in this case,
$C_{F,-c}$
behaves well, as the next Proposition~shows.
\end{iii}
\end{rems}

\begin{con}
We~adopt here the usual convention that
$\Psi$
is determined only up to inner~automorphisms.
\end{con}

\begin{rem}
An~experiment shows that for some but not all~subgroups
$$g \subseteq (S_2 \wr S_6) \cap A_{12} \cong U_{63} \subset G \, ,$$
the groups
$g$
and
$\Psi(g)$
are conjugate
in~$G$.
\end{rem}

\begin{prop}
\label{conj}
Let\/~$k$
be a field and\/
$F(T) = T^6 + a_5T^5 + \cdots + a_1T + a_0 \in k[T]$
be a separable, monic polynomial of degree six, where\/
$a_0 = c^2$,
for
some\/~$c \in k \setminus \{0\}$.

\begin{abc}
\item
Then the two conic~bundles
$$\pro\colon B_{F,c} \to \Pb^1 \qquad\text{and}\qquad \pro\colon B_{F,-c} \to \Pb^1$$
both split over the splitting
field\/~$l$
of\/~$F(T^2)$.
\item
The~operations of\/
$\Gal(l/k)$
on the components of the reducible fibres, however, differ by the outer automorphism\/~$\Psi$
of\/~$(S_2 \wr S_6) \cap A_{12}$.
\end{abc}\smallskip

\noindent
{\bf Proof.}
{\em
a)
is clear from Proposition~\ref{cb_gal}.b).ii).\smallskip

\noindent
b)
{\em First step.}
The relative situation.

\noindent
The construction
of~$B_{F,c}$,
as described in Proposition~\ref{cb_gal}, assigns to a monic polynomial
$F(T) = T^6 + a_5T^5 + \cdots + a_0$
and a choice of
$c$
such that
$c^2 = a_0$
a conic bundle with exactly six singular~fibres. As~noticed in Remark~\ref{discr}.ii), the whole process may be carried out in a relative situation, such that the old construction reappears when working fibre-by-fibre.

Concretely,~this means the~following. The~affine
space~$\Ab^6_k$
with coordinate functions
$a_0, \ldots, a_5$
admits the ramified double cover, given by
$w^2 = a_0$,
which is again an affine space
$\Spec k[w, a_1, \ldots, a_5] \cong \Ab^6_k$.
Over~a certain open subscheme
$\bm{W} \subset \Spec k[w, a_1, \ldots, a_5]$
to be specified below, the construction provides a family of conic bundles, which is the hypersurface
$$\bm{B} \subset \bm{W} \times \Pb^1 \times \Pb^2 \,,$$
given by
$(T_0\, T_1\, T_2)\, M \,(T_0\, T_1\, T_2)^t = 0$.
Here,~$M$
is the
$3\times3$-matrix~(\ref{matrix}).

The equation
$\det M = 0$
defines the locus
$\bm{L} \subset \bm{W} \!\times\! \Pb^1$,
over which the conics degenerate to rank two, i.e.~to the union of two~lines. Degeneration to even lower ranks does not~occur. The~locus
$\bm{L}$
is, in fact, contained in~$\bm{W} \!\times\! \Ab^1$
and given by the equation
$T^6 + a_5T^5 + \cdots + a_1T + w^2 = 0$.

The~lines themselves hence form a
\mbox{$\Pb^1$-bundle}
over a double cover
${\bm X}$
of~$\bm{L}$.
Over~$\bm{L}$,
${\bm X}$~is
given by the equation
$W^2 = \det M_{33}$,
for~$M_{33}$
the upper left
$2\times2$-minor
of~$M$.
It~turns out that
$\det M_{33}$
coincides
with~$T$,
up to a factor being a perfect~square.\smallskip

\noindent
{\em Second step.}
The base scheme.

\noindent
In~the affine
space~$\Ab^6_k$
with coordinate functions
$a_0, \ldots, a_5$,
there is the discriminant locus
$\Delta \subset \Ab^6_k$
of the polynomial
$T^{12} + a_5T^{10} + \cdots + a_0$.
This~is a reducible divisor consisting of two~components. One~is the hyperplane
$Z(a_0)$,
the other is the discriminant locus of the~polynomial
$T^6 + a_5T^5 + \cdots + a_0$.

We~specify
$\bm{W}$
to be the preimage
of~$\Ab^6_k \setminus \Delta$
under the double
cover~$w = a_0^2$.
In~fact, in Proposition~\ref{cb_gal}, we only worked with points on the somewhat smaller base scheme
$\Ab^6_k \setminus \Delta \setminus \Delta_\sing$,
but this will not make any difference for the argument~below.

There~are the finite \'etale morphisms, i.e.\ \'etale covers,
$$\bm{X} \stackrel{b}{\longrightarrow} \bm{W} \stackrel{q}{\longrightarrow} \Ab^6_k \setminus \Delta$$
of
\mbox{$k$-schemes}.
Here,
$q\colon \bm{W} \to \Ab^6_k \setminus \Delta$
is the structural double cover~morphism.
Moreover,~$b\colon \bm{X} \to \bm{W}$
is the composition of the natural projection
$\bm{X} \twoheadrightarrow \bm{L}$
with the composition
$$\bm{L} \hookrightarrow \bm{W} \!\times\! \Pb^1 \stackrel{\pr_1}{\twoheadrightarrow} \bm{W}.$$
According~to its construction,
$\bm{X}$~is
isomorphic to the integral closure \cite[Section~6.3]{EGAII} of
$\bm{W}$
relative to the extension
$Q(\bm{W})[T]/(T^{12} + a_5T^{10} + \cdots + a_1T^2 + w^2)$
of the function
field~$Q(\bm{W})$.
Moreover,~under this isomorphism, 
$b\colon \bm{X} \to \bm{W}$
goes over into the structural~morphism. Since~the discriminant locus has been taken out,
$b\colon \bm{X} \to \bm{W}$
is indeed \'etale~\cite[Exp.~I, Corollaire~7.4]{SGA1}.

The~twelve lines occurring in the conic bundle
$B_{F,c}$
are in a natural bijection with the twelve geometric points on the corresponding fibre
of~$b$.
Hence,~in what follows it suffices to consider the \'etale cover
$b\colon \bm{X} \to \bm{W}$
instead of the lines~themselves.\smallskip

\noindent
{\em Third step.}
The generic point.

\noindent
Over the generic point
of~$\Ab^6_k$
with function field
$\bm{F} = k(a_0, \ldots, a_5)$
lies exactly one point
of~$\bm{W}$,
the generic point
$\eta$,
corresponding to the function~field
$$\bm{F}(\sqrt{\mathstrut a_0}) = k(\sqrt{\mathstrut a_0}, a_1, \ldots, a_5) = k(w, a_1, \ldots, a_5) \, .$$
Above~$\eta$,
there is the generic fibre
$\bm{X}_\eta$,
which is again only one point and corresponds to the function~field
$$Q(\bm{X}_\eta) = k(w, a_1, \ldots, a_5)[T]/(T^{12} + a_5T^{10} + \cdots + a_1T^2 + w^2) \, .$$

The Galois group of the polynomial
$T^{12} + a_5T^{10} + \cdots + a_1T^2 + a_0$
over~$\bm{F}$
is equal to
$S_2 \wr S_6$,
as the coefficients
$a_0, \ldots, a_5$
are indeterminates. Moreover, the discriminant
of~$T^{12} + a_5 T^{10} + \cdots + a_0$
is~$a_0$,
up to factors being squares
in~$\bm{F}$.
Thus,
$$\bm{F}(\sqrt{\mathstrut a_0}) \subset \Split_{\bm{F}} (T^{12} + a_5 T^{10} + \cdots + a_1T^2 + a_0) \, ,$$
which reduces the Galois group
over~$\smash{\bm{F}(\sqrt{\mathstrut a_0})}$
to only even~permutations.

As~a consequence of this, the Galois hull
$\overline{\bm{X}}$
of~$\bm{X}$
over
$\bm{W}$
is automatically Galois
over~$\Ab^6_k \setminus \Delta$.
The~\'etale~cover
$$q \!\circ\! \overline{b}\colon \overline{\bm{X}} \stackrel{\overline{b}}{\longrightarrow} \bm{W} \stackrel{q}{\longrightarrow} \Ab^6_k \setminus \Delta$$
of Galois group
$\Aut(\overline{\bm{X}}) \cong S_2 \wr S_6$
is split into two parts, both of which are \'etale covers and Galois, with Galois groups
$\Aut({\bm{W}}) \cong \bbZ/2\bbZ$
for~$q$
and
$\Aut_{\bm{W}}(\overline{\bm{X}}) \cong (S_2 \wr S_6) \cap A_{12}$
for~$\smash{\overline{b}}$.\smallskip

\noindent
{\em Fourth step.}
\smash{\'Etale} fundamental groups.

\noindent
Now let
$F \in k[T]$
be a monic polynomial of degree six that is separable and has constant
term~$a_0 = c^2$.
The~coefficients of this polynomial define two
\mbox{$k$-rational}
points
$\smash{(F,c)}$,
$\smash{(F,-c) \in \bm{W}(k)}$
lying over the same point\vspace{.2mm}
$\smash{(F,c^2) \in (\Ab^6_k \setminus \Delta)(k)}$.

Let~us fix an algebraic closure
$\smash{\overline{k}}$
and geometric points
$\smash{\overline{(F,c)}\colon \Spec \overline{k} \to \bm{W}}$
and
$\smash{\overline{(F,-c)}\colon \Spec \overline{k} \to \bm{W}}$,
lying above
$(F,c)$
and
$(F,-c)$,
respectively, that are mapped
under~$q$
to the same geometric point
$\smash{\overline{(F,c^2)}}$
on~$\Ab^6_k \setminus \Delta$.
We~also fix two geometric points
$p_c$
and~$p_{-c}$
on~$\smash{\overline{\bm{X}}}$,
which are mapped
under~$\smash{\overline{b}}$
to
$\smash{\overline{(F,c)}}$
and~$\smash{\overline{(F,-c)}}$,
respectively, and shall serve as base points for the operations of the \'etale fundamental groups on the \'etale covers
$\smash{\overline{b}}$
and~$\smash{q \!\circ\! \overline{b}}$.

Furthermore,~there are the natural homomorphisms
\begin{eqnarray*}
q_*^{(c)} \colon \pi_1^\et(\bm{W}, \overline{(F,c)}) & \longrightarrow & \pi_1^\et(\Ab^6_k \setminus \Delta, \overline{(F,c^2)}) \quad\text{and} \\
q_*^{(-c)} \colon \pi_1^\et(\bm{W}, \overline{(F,-c)}) & \longrightarrow & \pi_1^\et(\Ab^6_k \setminus \Delta, \overline{(F,c^2)})
\end{eqnarray*}
between the \'etale fundamental groups~\cite[Exp.~V]{SGA1}, which are injective and onto a normal subgroup of
index~$2$.
Indeed,~$q$
is an \'etale double cover and
$\bm{W}$
is~connected.

In~addition, the \'etale cover
$q \!\circ\! \overline{b}$
corresponds, under the equivalence of categories described in~\cite[Exp.~V, Section 7]{SGA1}, to a surjective continuous group homomorphism
$$\varrho\colon \pi_1^\et(\Ab^6_k \setminus \Delta, \overline{(F,c^2)}) \longrightarrow \Aut(\overline{\bm{X}}) \, .$$
Similarly,~the \'etale cover
$\overline{b}$
yields the surjective continuous homomorphisms
$$\varrho_c\colon \pi_1^\et(\bm{W}, \overline{(F,c)}) \to \Aut_{\bm{W}}(\overline{\bm{X}}) \;\;\text{ and }\;\; \varrho_{-c}\colon \pi_1^\et(\bm{W}, \overline{(F,-c)}) \to \Aut_{\bm{W}}(\overline{\bm{X}}) \, ,$$
respectively. Both~are compatible
with~$\varrho$
in the sense that the diagrams
$$
\xymatrix{
\pi_1^\et(\bm{W}, \overline{(F,c)}) \ar@{->>}[r]^{\;\;\;\;\varrho_c} \ar@{^{(}->}[d]_{q_*^{(c)}} & \Aut_{\bm{W}}(\overline{\bm{X}}) \ar@{^{(}->}[d]^{\text{incl.}} \\ 
\pi_1^\et(\Ab^6_k \setminus \Delta, \overline{(F,c^2)}) \ar@{->>}[r]^{\;\;\;\;\;\;\;\;\;\varrho} & \Aut(\overline{\bm{X}})
}
\quad\raisebox{-8mm}{\text{ and }}\quad
\xymatrix{
\pi_1^\et(\bm{W}, \overline{(F,-c)}) \ar@{->>}[r]^{\;\;\;\;\varrho_{-c}} \ar@{^{(}->}[d]_{q_*^{(-c)}} & \Aut_{\bm{W}}(\overline{\bm{X}}) \ar@{^{(}->}[d]^{\text{incl.}} \\ 
\pi_1^\et(\Ab^6_k \setminus \Delta, \overline{(F,c^2)}) \ar@{->>}[r]^{\;\;\;\;\;\;\;\;\;\varrho} & \Aut(\overline{\bm{X}})
}
$$
commute.

Next,~we have to match the two operations
$\varrho_c$
and~$\varrho_{-c}$
against each~other. For~this, we choose an isomorphism
$\smash{\iota\colon \Phi_{\overline{(F,-c)}} \to \Phi_{\overline{(F,c)}}}$
of fibre functors, i.e.\ a ``homotopy class of paths''
$$s\in \pi_1^\et(\bm{W}, \overline{(F,c)}, \overline{(F,-c)})$$
in the fundamental groupoid, cf.~\cite[Exp.~V, Section 7]{SGA1} or~\cite[Paragraph~10.16]{De}.
We~assume that
$s$~lifts
to a path
in~$\pi_1^\et(\overline{\bm{X}}, p_c, p_{-c})$.
The~class~$s$
defines an~isomorphism
\begin{eqnarray*}
\iota_s\colon \pi_1^\et(\bm{W}, \overline{(F,c)}) & \longrightarrow & \pi_1^\et(\bm{W}, \overline{(F,-c)}) \,, \\
 \sigma & \mapsto & s \!\circ\! \sigma \!\circ\! s^{-1},
\end{eqnarray*}
by~conjugation.

Fortunately,
$q_*(s) \in \pi_1^\et(\Ab^6_k \setminus \Delta, \overline{(F,c^2)})$
is an element of an ordinary \'etale fundamental group. Thus,~the commutativity of the diagrams above yields~that
\begin{eqnarray*}
\varrho_{-c}(\iota_s(\sigma)) = \varrho_{-c}(s \!\circ\! \sigma \!\circ\! s^{-1}) = \varrho\big(q_*^{(-c)}(s \!\circ\! \sigma \!\circ\! s^{-1})\big) = \varrho\big(q_*(s) \!\circ\! q_*^{(-c)}(\sigma) \!\circ\! q_*(s)^{-1}\big) = \hspace{.3cm} \\
\varrho(q_*(s)) \!\cdot\! \varrho(q_*^{(c)}(\sigma)) \!\cdot\! \varrho(q_*(s))^{-1} = \varrho(q_*(s)) \!\cdot\! \varrho_c(\sigma) \!\cdot\! \varrho(q_*(s))^{-1},
\end{eqnarray*}
for every
$\sigma \in \pi_1^\et(\bm{W}, \overline{(F,c)})$.
I.e.,
$\varrho_{-c} \!\circ\! \iota_s$
and
$\varrho_c$
differ by conjugation
with~$\varrho(q_*(s))$,
$$\varrho_{-c} \!\circ\! \iota_s = \varrho(q_*(s)) \cdot\! \varrho_c \!\cdot \varrho(q_*(s))^{-1} \, .$$

We~note finally that, according to its construction,
$q_*(s)$
is an element of the fundamental group
$\smash{\pi_1^\et(\Ab^6_k \setminus \Delta, \overline{(F,c^2)})}$
that does not lift to the fundamental group
of~$\bm{W}$.
Hence,~it lies in the only nontrivial coset
$\smash{\pi_1^\et(\Ab^6_k \setminus \Delta, \overline{(F,c^2)}) \,\setminus\, q_*^{(c)}\big( \pi_1^\et(\bm{W}, \overline{(F,c)}) \big)}$
and one~has
$$\varrho(q_*(s)) \in \Aut(\overline{\bm{X}}) \setminus \Aut_{\bm{W}}(\overline{\bm{X}}) \,.$$
In~other words, the operations
$\varrho_{-c} \!\circ\! \iota_s$
and
$\varrho_c$
indeed differ by the outer automorphism
$\Psi$
of
$\Aut_{\bm{W}}(\overline{\bm{X}}) \cong (S_2 \wr S_6) \cap A_{12}$.\smallskip\pagebreak[3]

\noindent
{\em Fifth step.}
Conclusion--The absolute Galois group.

\noindent
In~order to complete the proof, we still have to specify which class of paths
$\smash{s\in \pi_1^\et(\bm{W}, \overline{(F,c)}, \overline{(F,-c)})}$
to work with. For~this, let us write
$k^\sep$
for the separable closure
of~$k$
within the chosen algebraic
closure~$\overline{k}$.
Then~the
\mbox{$k$-rational}
point
$(F,-c)\colon \Spec k \to \bm{W}$
yields a homomorphism
$$i_{-c}\colon \Gal(k^\sep/k) = \pi_1^\et(\Spec k, \Spec \overline{k}) \longrightarrow \pi_1^\et(\bm{W}, \overline{(F,-c)}) \, .$$
Similarly,~the
\mbox{$k$-rational}
point
$(F,c)\colon \Spec k \to \bm{W}$
provides a homomorphism
$\smash{i_c\colon \Gal(k^\sep/k) = \pi_1^\et(\Spec k, \Spec \overline{k}) \longrightarrow \pi_1^\et(\bm{W}, \overline{(F,c)})}$.

The compositions
$\varrho_{-c} \!\circ\!i_{-c}$
and~$\varrho_c \!\circ\! i_c$
then describe the operations
of~$\Gal(k^\sep/k)$
on~the \'etale
cover~$\smash{\overline{b}}$,
obtained by taking
$p_{-c}$
and~$p_c$,
respectively, as base~points. I.e., the Galois operation on the fibres
$\smash{\overline{b}^{-1}(\overline{(F,-c)})}$
and~$\smash{\overline{b}^{-1}(\overline{(F,c)})}$.
Both these fibres are transitive
$\Gal(k^\sep/k)$-sets
isomorphic to
$$\Gal(k^\sep/k)/\Gal(k^\sep/\Split_k (F(T^2))) \, .$$
Taking~$p_{-c}$
and~$p_c$
as generators, we thus find a bijection
$\smash{\overline{b}^{-1}(\overline{(F,-c)}) \cong \overline{b}^{-1}(\overline{(F,c)})}$
that is compatible with the Galois~operations. This~may be extended to an isomorphism between the fibre functors, i.e.\ to a class
$\smash{s\in \pi_1^\et(\bm{W}, \overline{(F,c)}, \overline{(F,-c)})}$,
as~required.

By~construction, the homomorphisms
$\smash{i_{-c}, \iota_s \!\circ\! i_c\colon \Gal(k^\sep/k) \to \pi_1^\et(\bm{W}, \overline{(F,-c)})}$
agree modulo elements of
$\smash{\pi_1^\et(\bm{W}, \overline{(F,-c)})}$
acting trivially
on~$\smash{\overline{\bm{X}}}$.
Thus,~the result of the previous step implies~that
$$\varrho_{-c} \!\circ\! \iota_s \!\circ\! i_c = \varrho_{-c} \!\circ\! i_{-c} \colon \Gal(k^\sep/k) \to \Aut_{\bm{W}}(\overline{\bm{X}})
\text{ and }
\varrho_c \!\circ\! i_c \colon \Gal(k^\sep/k) \to \Aut_{\bm{W}}(\overline{\bm{X}})$$
differ by the outer
automorphism~$\Psi$
of
$\smash{\Aut_{\bm{W}}(\overline{\bm{X}}) \cong (S_2 \wr S_6) \cap A_{12}}$,
as~claimed.
}
\eop
\end{prop}

\begin{theo}
\label{main_cb}
Let an infinite
field\/~$k$
of characteristic
not\/~$2$,
a normal and separable extension
field\/~$l$,
and an injective group~homomorphism
$$i\colon \Gal(l/k) \hookrightarrow U_{63} \subset G \cong \Sp_6(\bbF_{\!2}) \subset S_{28}$$
be given. Then there exists a nonsingular quartic
curve\/~$C$
over\/~$k$
such that\/
$l$
is the field of definition of the 28 bitangents and each\/
$\sigma \in \Gal(l/k)$
permutes the bitangents as described by\/
$i(\sigma) \in G \subset S_{28}$.\medskip

\noindent
{\bf Proof.}
{\em
There is an isomorphism
$U_{63} \cong (S_2 \wr S_6) \cap A_{12}$,
which describes the operation
of~$U_{63}$
on the twelve lines of the distinguished Steiner~hexad.

We~are thus given an injection
$\smash{i'\colon \Gal(l/k) \hookrightarrow (S_2 \wr S_6) \cap A_{12}}$,
to which Lemma \ref{pol_12} can be applied. It~yields a polynomial
$F = T^6 \!+\! a_5T^5 \!+ \cdots +\! a_1T \!+\! a_0 \in k[T]$
of degree six, such that the splitting field
of~$F(T^2)$
is
exactly~$l$.
More~explicitly,
$\smash{l = k(\sqrt{A_1}, ..., \sqrt{A_6})}$,
for~$A_1, \ldots, A_6$
the roots
of~$F$.
Furthermore,~the operation of
$\Gal(l/k)$
on the square roots
$\pm\sqrt{A_i}$
agrees with the natural operation of
$S_2 \wr S_6$
on twelve objects forming six~pairs. As~discussed in Corollary~\ref{Zariski}, there is a Zariski dense subset
of~$k[T]_6$
of polynomials with the same~behaviour.

As~$\smash{\Gal(l/k) \stackrel{i'}{\hookrightarrow} A_{12}}$,
we know that the discriminant
of~$F(T^2)$
is a square
in~$k$.
In~other words,
$\smash{a_0 \in (k^*)^2}$.
Hence,~we may apply Proposition~\ref{cb_gal} in order to find two nonsingular plane quartics
$C_{F,c}$
and
$C_{F,-c}$.
According~to Lemma~\ref{crit_hex}, they both enjoy the following~property. Twelve~of their bitangents, which form Steiner hexads
$\Hb^+$
and~$\Hb^-$,
respectively, are defined
over~$l$
and permuted by
$\Gal(l/k)$
exactly in the way described by the given injection
$\smash{\Gal(l/k) \stackrel{i'}{\hookrightarrow} S_2 \wr S_6}$.

Let~us consider
$C_{F,c}$~first.
Having~equipped the 28 bitangents
of~$C_{F,c}$
with any numbering distinguishing the Steiner
hexad~$\Hb^+$,
the operation
of~$\Gal(l/k)$
on~$\Hb^+$
agrees with the desired one up to an inner automorphism
$\varphi$
of~$S_2 \wr S_6$.
This~might be the conjugation with an even element. In~this case,
$\varphi$~extends
to an inner automorphism of the whole
of~$G$,
which completes the argument, as discussed in~Remark~\ref{mark}.i). Otherwise,~Proposition~\ref{conj} shows that
$C_{F,-c}$
has all the properties~required.
}
\eop
\end{theo}

\begin{rem}[The case that
$k$
is a number field]

Let~$k$
be a number field and
$g$
a subgroup of
$G \cong \Sp_6(\bbF_{\!2})$
that is contained
in~$U_{63}$.
Then there exists a nonsingular quartic
curve~$C$
over~$k$
such that the natural permutation~representation
$$i\colon \Gal(\overline{k}/k) \longrightarrow G \subset S_{28}$$
on the 28 bitangents
of~$C$
has the
subgroup~$g$
as its~image.\medskip

\noindent
{\bf Proof.}
According~to Theorem~\ref{main_cb}, it suffices to show that, for every number field~$k$
and each
subgroup~$g \subseteq U_{63}$,
there exists a normal extension
field~$l$
such that
$\Gal(l/k)$
is isomorphic
to~$g$.
This~is a particular instance of the inverse Galois problem, but the groups occurring are easy~enough.

In~fact, among the 1155 conjugacy classes of subgroups
of~$G$
that are contained
in~$U_{63}$,
1119 consist of solvable~groups. For~these, the inverse Galois problem has been solved by I.\,R.~Shafarevich~\cite{Sha}, cf.~\cite[Theorem~9.5.1]{NSW}.
The~groups in the remaining 36 conjugacy classes all turn out to be factor groups of the wreath product
$S_2 \wr H$,
where
$H \subset S_{15}$
is isomorphic to
$A_5$,
$S_5$,
$A_6$,
or~$S_6$,
cases for which Galois extensions are known over an arbitrary number field~\cite[paragraph~2.2.3]{KM}.%
\eop
\end{rem}

\begin{ex}
Put
$$F := T^6 - 24T^5 + 152T^4 - 340T^3 + 335T^2 - 150T + 25 \in \bbQ[T] \, .$$
Then the construction described in Proposition~\ref{cb_gal} yields the plane~quartics
\begin{align*}
-2T_0^3T_2 + 37T_0^2T_1T_2 + 67T_0^2T_2^2 + 2T_0T_1^3 - 10T_0T_1^2T_2 + 114T_0T_1T_2^2 + 166T_0T_2^3 + 4T_1^4 \\
{} - 168T_1^3T_2 + 42T_1^2T_2^2 + 369T_1T_2^3 - 45T_2^4 = 0
\end{align*}
for
$c = 5$~and
\begin{align*}
-2T_0^3T_2 + 43T_0^2T_1T_2 + 480T_0^2T_2^2 + 2T_0T_1^3 - 36T_0T_1^2T_2 - 
2460T_0T_1T_2^2 + 8700T_0T_2^3 \\
{} + 2T_1^4 - 189T_1^3T_2 + 5170T_1^2T_2^2 - 31255T_1T_2^3 + 40250T_2^4 = 0
\end{align*}
for~$c = -5$.
In~this case, the two subgroups
of~$G$
occurring as the Galois groups operating on the 28~bitangents are in fact conjugate to each~other.

The Galois group of the polynomial
$F$
itself is
$A_4$,
realised as a subgroup
of~$S_6$
by the operation on 2-sets. On~the other hand, the Galois group of
$F(T^2)$
is of
index~$8$
in
$(S_2 \wr A_4) \cap A_{12}$,
of
order~$48$.
According~to the classification of transitive groups in degree twelve, due to G.\,F.\ Royle~\cite{Ro} (cf.~\cite{CHM}) and used by {\tt magma} as well as {\tt gap}, it corresponds to number~12T31.
\end{ex}

\begin{ex}
Over~the function field
$\bbF_{\!3}(t)$,
the~polynomial
\begin{align*}
T^6 + (2t^4 + 2t^2)T^4 + (2t^6 + 2t^4 + t^2 + 1)&T^3 + (t^8 + t^6 + 2t^4 +
    2t^2)T^2 \\[-1mm]
{} &+ (2t^{10} + 2t^4)T + t^{12} + 2t^6 + 1 \in \bbF_{\!3}(t)[T]
\end{align*}
provides the same Galois~group. We~obtain the nonsingular plane quartics over
$\bbF_{\!3}(t)$,
given~by
\begin{align*}
& T_0^4 + 2T_0^3T_1 + (2t^6 \!+\! t^4 \!+\! t^2)T_0^3T_2 + (t^4 \!+\! t^2 \!+\! 1)T_0^2T_1T_2 + (2t^4 \!+\! t^2 \!+\! 2)T_0^2T_2^2 \\[-1mm]
& {} + T_0T_1^2T_2 + (2t^6 \!+\! 2t^4 \!+\! 2)T_0T_1T_2^2 + (t^{10} \!+\! t^8 \!+\! 2t^4 \!+\! 2)T_0T_2^3 + T_1^3T_2 \\[-1mm]
& {} \hspace{5cm} + (2t^4 \!+\! 2)T_1^2T_2^2 + (t^8 \!+\! t^6 \!+\! t^2 \!+\! 1)T_1T_2^3 + 2t^2T_2^4 = 0
\end{align*}
and
\begin{align*}
& T_0^4 + 2T_0^3T_1 + (t^6 \!+\! t^4 \!+\! t^2 \!+\! 2)T_0^3T_2 + (t^4 \!+\! t^2 \!+\! 1)T_0^2T_1T_2 \\[-1mm]
& {} + (2t^{10} \!+\! 2t^8 \!+\! 2t^6 \!+\! t^4 \!+\! 1)T_0^2T_2^2 + T_0T_1^2T_2 + (2t^6 \!+\! 2t^4 \!+\! 2)T_0T_1T_2^2 \\[-1mm]
& {} + (t^{12} \!+\! 2t^{10} \!+\! 2t^2 \!+\! 1)T_0T_2^3 + T_1^3T_2 + (2t^6 \!+\! 2t^4 \!+\! 1)T_1^2T_2^2 \\[-1mm]
& {} + (t^{10} \!+\! 2t^8 \!+\! 2t^6 \!+\! t^4 \!+\! 2t^2 \!+\! 2)T_1T_2^3 + (2t^{16} \!+\! 2t^{12} \!+\! t^{10} \!+\! t^6 \!+\! 2t^4 \!+\! 2t^2 \!+\! 2)T_2^4 = 0 \, .
\end{align*}
\end{ex}

\begin{ex}
Put
$$F := T^6 - 3T^5  -2T^4  + 9T^3 - 5T^1 + 1 \in \bbQ[T] \, .$$
Then the construction described in Proposition~\ref{cb_gal} yields the plane~quartics
\begin{align*}
-T_0^3T_2 - 2T_0^2T_1T_2 + 14T_0^2T_2^2 - 2T_0T_1^3 + 9T_0T_1^2T_2 + 4T_0T_1T_2^2 - 7T_0T_2^3 + T_1^4 - 2T_1^3T_2 \\
{} - 7T_1^2T_2^2 + 3T_2^4 = 0
\end{align*}
for
$c = 1$~and
\begin{align*}
-T_0^3T_2 - 2T_0^2T_1T_2 + 6T_0^2T_2^2 - 2T_0T_1^3 + 9T_0T_1^2T_2 + 12T_0T_1T_2^2 - 3T_0T_2^3 + T_1^4 + 6T_1^3T_2 \\
{} + 5T_1^2T_2^2 - 4T_1T_2^3 - T_2^4 = 0
\end{align*}
for~$c = -1$.
In~this case, the two subgroups
of~$G$
occurring as the Galois groups on the 28~bitangents are not conjugate to each~other. For~example, the second is contained in the subgroup
$U_{36} \subset G$,
too, while the first one is~not.

The Galois group of the polynomial
$F$
itself is
$S_3$,
realised as a subgroup
of~$S_6$
by the regular~representation. The~Galois group of
$F(T^2)$
is of
index~$2$
in
$(S_2 \wr S_3) \cap A_{12}$,
of
order~$96$.
According~to the classification of transitive groups in degree twelve, it is number~12T69.
\end{ex}

\begin{ex}
Over~the function field
$\bbF_{\!3}(t)$,
the~polynomial
$$T^6 + 2t^2T^5 + t^4T^4 + (2t^7 + t^4 + t^2 + 1)(T^3 + t^2T^2) + (t^7 + t^4 +
t^2 + 1)^2 \in \bbF_{\!3}(t)[T]$$
provides the same Galois~group. We~obtain the nonsingular plane quartics over
$\bbF_{\!3}(t)$,
given~by
\begin{align*}
& T_0^4 + 2T_0^3T_1 + (2t^7 \!+\! 2t^4)T_0^3T_2 + (t^2 \!+\! 1)T_0^2T_1T_2 + (t^9 \!+\! 2t^7 \!+\! 2t^4 \!+\! 2t^2 \!+\! 2)T_0^2T_2^2 \\[-1mm]
& {} + (2t^2 \!+\! 1)T_0T_1^2T_2 + (2t^6 \!+\! t^4 \!+\! 2)T_0T_1T_2^2 + (2t^{11} \!+\! 2t^9 \!+\! t^7 \!+\! t^6 \!+\! t^4 \!+\! t^2 \!+\! 2)T_0T_2^3 \\[-1mm]
& {} + T_1^3T_2+ (t^2 \!+\! 2)T_1^2T_2^2 + (t^9 \!+\! t^7 \!+\! t^6 \!+\! 2t^4 \!+\! 2t^2 \!+\! 1)T_1T_2^3 \\[-1.1mm]
& {} \hspace{6.5cm} + (t^{14} \!+\! 2t^{13} \!+\! t^{12} \!+\! t^{11} \!+\! 2t^9 \!+\! 2t^6 \!+\! t^4)T_2^4 = 0
\end{align*}
and
\begin{align*}
& T_0^4 + 2T_0^3T_1 + (t^7 \!+\! t^4 \!+\! 2t^2 \!+\! 2)T_0^3T_2 + (t^2 \!+\! 1)T_0^2T_1T_2 + (t^7 \!+\! 2t^6 \!+\! 1)T_0^2T_2^2 \\[-1mm]
& {} + (2t^2 \!+\! 1)T_0T_1^2T_2 + (t^9 \!+\! 2t^4 \!+\! t^2 \!+\! 2)T_0T_1T_2^2 \\[-1mm]
& {} + (t^{14} \!+\! 2t^{11} \!+\! 2t^9 \!+\! 2t^8 \!+\! t^7 \!+\! 2t^6 \!+\! t^4 \!+\! 2t^2 \!+\! 1)T_0T_2^3 + T_1^3T_2 + (2t^7 \!+\! 2t^4 \!+\! 1)T_1^2T_2^2 \\[-1mm]
& {} + (2t^9 \!+\! 2t^7 \!+\! 2t^6 \!+\! t^4 \!+\! t^2 \!+\! 2)T_1T_2^3 \\[-1mm]
& {} \hspace{4.5cm} + (t^{18} \!+\! 2t^{14} \!+\! 2t^{13} \!+\! t^{11} \!+\! t^{10} \!+\! 2t^9 \!+\! 2t^8 \!+\! t^6 \!+\! t^2 \!+\! 2)T_2^4 = 0 \, .
\end{align*}
\end{ex}

\begin{rem}
The examples above were chosen from the enormous supply in the hope that they are of some particular~interest. The Galois groups realised are in fact the minimal ones that yield the generic orbit type
$[12,16]$
on the 28~bitangents. The~corresponding conic bundle surfaces are of the minimal possible Picard
rank~$2$,
and their generic quadratic twists are of Picard
rank~$1$.
\end{rem}

\section{Twisting}
\label{fuenf}

There~is the double cover
$\smash{p\colon \widetilde{G} \to G}$
of finite groups, for
$\smash{W(E_7) \cong \widetilde{G} \subset S_{56}}$
and
$\Sp_6(\bbF_{\!2}) \cong G \subset S_{28}$,
which is given by the operation on the size two~blocks. The~kernel
of~$p$
is exactly the
centre~$\smash{Z \subset \widetilde{G}}$.
For~a subgroup
$\smash{H \subset \widetilde{G}}$,
one therefore has two~options.

\begin{iii}
\item
Either~$p|_H \colon H \to p(H)$
is two-to-one. Then
$H = p^{-1}(h)$,
for~$h := p(H)$.
In~this case,
$H$~contains
the centre
of~$\smash{\widetilde{G}}$
and, as abstract groups, one has an isomorphism
$H \cong p(H) \times \bbZ/2\bbZ$.
\item
Or~$p|_H \colon H \to p(H)$
is bijective.
\end{iii}
In~our geometric setting, the first case is the generic~one. More~precisely, let
$C\colon q=0$
be a nonsingular plane quartic such that the 28~bitangents are acted upon by the
group~$h \subseteq G$.
Then,~for
$\lambda$
an indeterminate, the 56~exceptional curves on
$S_\lambda\colon \lambda w^2 = q$
operated upon
by~$h \times \bbZ/2\bbZ$.

\begin{lem}
\label{univ}
Let\/~$k$
be any field and\/
$C\colon q = 0$
a nonsingular plane quartic
over\/~$k$.
Write\/~$l_0$
for the field of definition of the 28~bitangents
of\/~$C$
and let\/~$h \subset G$
be the subgroup, via which\/
$\smash{\Gal(k^\sep/k)}$
operates on~them.\smallskip

\noindent
Furthermore,~let
$S_t\colon t w^2 = q$
be the universal double cover
of\/~$\smash{\Pb^2_{k(t)}}$,
ramified
at\/~$C$,
defined over the function
field\/~$k(t)$.
Then~the following~holds.

\begin{abc}
\item
The~Galois
group\/~$\smash{\Gal(k(t)^\sep/k(t))}$
operates on the 56 exceptional curves
of\/~$S_t$
via\/~$p^{-1}(h)$.
\item
The~field of definition of the 56 exceptional curves
of\/~$S_t$
is\/
$\bm{l} = l_0(t)(\sqrt{ct})$,
for some\/
$c \in k^*$.
\end{abc}\smallskip

\noindent
{\bf Proof.}
{\em
a)
One~only has to show that the
field~$\bm{l}$
of definition of the 56 exceptional curves is strictly larger
than the
composite~$l_0 \cdot k(t) = l_0(t)$.
For~this, let us recall the local formula~(\ref{split_loc}), which shows that, indeed,
$\smash{\bm{l} = l_0(t)(\sqrt{ct})}$,
for a
certain~$c \in l_0^*$.\smallskip

\noindent
b)
The~field
$\bm{l}$
is necessarily Galois over
$k(t)$
and, as we are in case~i), the natural exact~sequence
$$0 \longrightarrow \Gal(\bm{l}/l_0(t)) \longrightarrow \Gal(\bm{l}/k(t)) \stackrel{\res}{\longrightarrow} \Gal(l_0(t)/k(t)) \longrightarrow 0$$
is split. I.e.,
$$\Gal(\bm{l}/k(t)) \cong \Gal(l_0(t)/k(t)) \times \bbZ/2\bbZ = \Gal(l_0/k) \times \bbZ/2\bbZ \, .$$
This~shows that
$\smash{\bm{l} = l_0(t)(\sqrt{p(t)})}$,
for some polynomial
$p(t) \in k[t]$.
On~the other hand, we just found that
$\smash{\bm{l} = l_0(t)(\sqrt{ct})}$,
for a certain constant
$c \in l_0^*$.
Both~results may be true, simultaneously, only if
$\smash{\bm{l} = l_0(t)(\sqrt{ct})}$,
where~$c$
is
in~$k^*$.
}
\eop
\end{lem}

For~particular choices
of~$\lambda$,
every subgroup
of~$\smash{\widetilde{G}}$
may be realised that has
image~$h$
under the
projection~$p$.

\begin{theo}
\label{twist}
Let~a
field\/~$k$
of characteristic
not\/~$2$,
a normal and separable extension
field\/~$l$,
and an injective group~homomorphism\/
$\smash{i\colon \Gal(l/k) \hookrightarrow \widetilde{G}}$
be~given.
Write\/~$l_0$
for the subfield corresponding to\/
$i^{-1}(Z)$
under the Galois correspondence.

\begin{abc}
\item
Then there is a commutative diagram
$$
\xymatrix{
\Gal(l/k) \ar@{^{(}->}[rr]^{\;\;\;\;\;\;\;\;i} \ar@{^{}->>}[d]_{\res} && \widetilde{G} \ar@{^{}->>}[d]^p \\
\Gal(l_0/k) \ar@{^{(}->}[rr]^{\;\;\;\;\;\;\;\;\overline\imath} && G \hsmash{ \, ,}
}
$$
the downward arrow on the left being the~restriction.
\item
Let\/~$C\colon q=0$
be a nonsingular plane quartic
over\/~$k$,
the 28 bitangents of which are defined
over\/~$l_0$
and acted upon
by\/~$\Gal(l_0/k)$
as described
by\/~$\overline\imath$.
Then~there exists some\/
$\lambda \in k^*$
such that the 56 exceptional curves of the degree two del Pezzo~surface
$$S_\lambda\colon \lambda w^2 = q$$
are defined
over\/~$l$
and each automorphism\/
$\sigma \in \Gal(l/k)$
permutes them as described by\/
$\smash{i(\sigma) \in \widetilde{G} \subset S_{56}}$.
\end{abc}\smallskip

\noindent
{\bf Proof.}
{\em
a)
follows directly from Galois~theory.\smallskip

\noindent
b)
Consider the universal double cover
$S_t$
over~$k(t)$,
given by
$tw^2 = q$.
According~to Lemma~\ref{univ}, this belongs to the first of the two cases distinguished~above. Moreover,~the field
$\bm{l}$
of definition of the 56 exceptional curves is
$\smash{\bm{l} = l_0(t)(\sqrt{ct})}$,
for
some~$c \in k^*$.
Let~us once again distinguish between the two cases, as~above.\smallskip\pagebreak[3]

\noindent
{\em First case.}
$\Gal(l/k) \cong \Gal(l_0/k) \times \bbZ/2\bbZ$.

\noindent
Then,~according to Kummer theory,
$\smash{l = l_0(\sqrt{\mathstrut a})}$
for some
$a \in k^*$.
Moreover~$\smash{\sqrt{\mathstrut a} \not\in l_0}$.
Specialising~$t$
to
$\lambda := ac$,
we find exactly the required field of definition.\smallskip

\noindent
{\em Second case.}
$l=l_0$.

\noindent
In~this case, a
subgroup~$H \subset \Gal(\bm{l}/k(t))$
of index two has been chosen that, under restriction, is mapped isomorphically
onto~$\Gal(l_0(t)/k(t))$.
We~want to
specialise~$t$
in such a way that we find
$H$
as the Galois group
over~$k$.

For this, we first observe that, under the Galois correspondence, the subgroup
$H$
corresponds to a quadratic extension field
$\bm{q}$
of~$k(t)$.
Thus,~one has
$H = \Gal(\bm{l}/\bm{q})$,
which yields that
$\bm{q} \not\subseteq l_0(t)$.
Indeed,~otherwise,
$H = \Gal(\bm{l}/\bm{q}) \supseteq \Gal(\bm{l}/l_0(t))$,
but
$\Gal(\bm{l}/l_0(t)) \cong \bbZ/2\bbZ$
is annihilated
under~$\res$.

On~the other hand, according to Galois~theory, the field extension
$\bm{l}/k(t)$,
i.e.\
$\smash{l_0(t)(\sqrt{ct})/k(t)}$,
has exactly three types of quadratic intermediate~fields. These~are

\begin{iii}
\item
the extension fields of the type
$\smash{k(t)(\sqrt{b})}$,
for
$\smash{k(\sqrt{b})}$
a quadratic subfield
of~$l_0$,
\item
the extension field
$\smash{k(t)(\sqrt{ct})}$,
\item
the extension fields of the type
$\smash{k(t)(\sqrt{b \!\cdot\! ct})}$,
for
$\smash{k(\sqrt{b})}$
a quadratic subfield
of~$l_0$.
\end{iii}

\noindent
Among~these, type~i) is excluded to us, as these fields are contained
in~$l_0(t)$.
Thus,~the subgroup
$H$
is realised as the Galois group over an intermediate field
$\smash{k(t)(\sqrt{b't})}$,
for some constant
$b' \in k^*$.
Specialising~$t$
to
$\lambda := b'$
yields the
subgroup~$H$
over~$k$.
}
\eop
\end{theo}

\begin{coro}
\label{main_dP2}
Let an infinite
field\/~$k$
of characteristic
not\/~$2$,
a normal and separable extension
field\/~$l$,
and an injective group~homomorphism
$$i\colon \Gal(l/k) \hookrightarrow p^{-1}(U_{63})$$
be~given. Then there exists a degree two del Pezzo
surface\/~$S$
over\/~$k$
such that\/
$l$
is the field of definition of the 56 exceptional curves and each\/
$\sigma \in \Gal(l/k)$
permutes them as described by\/
$\smash{i(\sigma) \in \widetilde{G} \subset S_{56}}$.\medskip

\noindent
{\bf Proof.}
{\em
This follows from Theorem~\ref{twist} together with Theorem~\ref{main_cb}.
}
\eop
\end{coro}

\begin{rem}
Let~$C\colon q^2 - q_1q_2 = 0$
(cf.~Proposition~\ref{detf}.a)) be a nonsingular plane quartic over a
field~$k$
with a Galois invariant Steiner hexad
and~$\lambda \in k^*$
a non-square.

\begin{iii}
\item
Then~the two conic bundles, associated with the Steiner hexad, on the twist
$S_\lambda\colon \lambda w^2 = q^2 - q_1q_2$
are no longer
\mbox{$k$-rational},
but defined over
$\smash{k(\sqrt{\lambda})}$
and conjugate to each~other.
\item
On~the other hand, in this case,
$S_\lambda$
carries a global Brauer~class
$\alpha \in \Br(S_\lambda)_2$,
which, over the function field
$k(S_\lambda)$,
is given by the quaternion algebra
$\smash{(\lambda, \frac{q_1}{T_0^2})_2}$,
for~$T_0$
one of the coordinate functions
on~$\Pb^2$.\medskip

\noindent
{\bf Proof.}
According~to \cite[Lemma~43.1.1 and Proposition~31.3]{Ma}, we only have to show that
$\div(\frac{q_1}{T_0^2})$
is the norm of a divisor
on~$(S_\lambda)_{k(\sqrt{\lambda})}$.
The equation
$$(q-\sqrt\lambda w)(q+\sqrt\lambda w) = q_1q_2$$
shows that
$Z(q-\sqrt\lambda w,q_1) - Z(T_0)$
indeed is such a divisor
on~$(S_\lambda)_{k(\sqrt{\lambda})}$.
\eop\medskip

\noindent
The 2-torsion Brauer classes on degree two del Pezzo surfaces have been systematically studied by P.~Corn in~\cite{Co}. It~turns out that the conjugacy classes of subgroups
of~$W(E_7)$
that lead to such a Brauer class form a partially ordered set with exactly two maximal elements. The~subgroup
$p^{-1}(U_{63})$
studied here is one of them. It yields the Brauer classes of the first type in Corn's~terminology.
\end{iii}
\end{rem}

\section{An application: Cubic surfaces with a Galois invariant double-six}
\label{sechs}

A~nonsingular cubic surface over an algebraically closed field contains exactly
27~lines. The~maximal subgroup
$\smash{G_\maxi \subset S_{27}}$
that respects the intersection pairing is isomorphic to the Weil
group~$W(E_6)$~\cite[Theorem~23.9]{Ma}
of
order~$51\,840$.

A~double-six (cf.~\cite[Remark~V.4.9.1]{Ha} or \cite[Subsection~9.1.1]{Do}) is a configuration of twelve lines
$E_1, \ldots, E_6,E'_1, \ldots, E'_6$
such that

\begin{iii}
\item
$E_1,\ldots,E_6$
are mutually skew,
\item
$E'_1,\ldots,E'_6$
are mutually skew, and
\item
$E_i \!\cdot\! E'_j = 1$,
for\/
$i \neq j$,
$1 \leq i,j \leq 6$,
and\/
$E_i \!\cdot\! E'_i = 0$
for~$i = 1, \ldots, 6$.
\end{iii}

\noindent
Every cubic surface contains exactly 36 double-sixes, which are transitively acted upon
by~$G_\maxi$
\cite[Theorem~9.1.3]{Do}. Thus,~there is an index-36 subgroup
$U_\ds \subset G_\maxi$
stabilising a double-six. This~is one of the maximal subgroups
of~$G_\maxi \cong W(E_6)$.
Up~to conjugation,
$W(E_6)$
has maximal subgroups of indices
$2$,
$27$,
$36$,
$40$,
$40$,
and~$45$.

As~an application of Theorem~\ref{twist}, we have the following~result.

\begin{theo}
\label{dsix}
Let an infinite
field\/~$k$
of characteristic
not\/~$2$,
a normal and separable extension
field\/~$l$,
and an injective group~homomorphism
$$i\colon \Gal(l/k) \hookrightarrow U_\ds \subset G_\maxi \cong W(E_6) \subset S_{27}$$
be~given. Then~there exists a nonsingular cubic
surface\/~$S$
over\/~$k$
such~that

\begin{iii}
\item
the 27 lines
on\/~$S$
are defined
over\/~$l$
and each\/
$\sigma \in \Gal(l/k)$
permutes them as described by\/
$i(\sigma) \in U_\ds \subset S_{27}$.
\item
$S$
is\/
\mbox{$k$-unirational}.
\end{iii}\smallskip

\noindent
{\bf Proof.}
{\em
There~is an injective homomorphism
$\iota\colon W(E_6) \hookrightarrow W(E_7)$
that corresponds to the blow-up of a point on a cubic~surface. The~subgroup
$\iota(U_\ds) \subset W(E_7)$
is of index
$36 \!\cdot\! 56 = 2016$.
It~is sufficient to show that the image
$\overline\iota(U_\ds)$
under the~composition
$$\overline\iota\colon W(E_6) \stackrel\iota\hookrightarrow W(E_7) \stackrel{p}\twoheadrightarrow W(E_7)/Z \cong G$$
is contained
in~$U_{63}$.
Indeed,~then Corollary~\ref{main_dP2} yields a degree two del Pezzo surface
$S'$
of degree two with a
\mbox{$k$-rational}
line~$L$
and the proposed Galois operation on the 27 lines that do not
meet~$L$.
Blowing~down
$L$,
we obtain a nonsingular cubic
surface~$S$
satisfying~i).
As~there is a
\mbox{$k$-rational}
point
on~$S$,
the blow down
of~$L$,
$S$~is
\mbox{$k$-unirational}~\cite[Theorem~2]{Ko02}.

To~show the group-theoretic claim, let us consider a model cubic surface over an algebraically closed~field. The~group
$U_\ds$
operates on the 27~lines with an orbit type
$[12, 15]$,
the orbit of size twelve being a double-six~\cite[Example~9.2]{EJ10}.
Hence,~$\iota(U_\ds)$
operates on the 56 exceptional curves of the degree two del Pezzo surface, obtained by blowing up one point, with orbit type
$[1,1,12,12,15,15]$.

A size twelve orbit consists of exceptional curves
$E_1, \ldots, E_6, E'_1, \ldots, E'_6$
such~that

\begin{iii}
\item
$E_1,\ldots,E_6$
are mutually skew,
\item
$E'_1,\ldots,E'_6$
are mutually skew, and
\item
$E_i \!\cdot\! E'_j = 1$,
for\/
$i \neq j$,
$1 \leq i,j \leq 6$,
and\/
$E_i \!\cdot\! E'_i = 0$
for~$i = 1, \ldots, 6$.
\end{iii}

\noindent
The~second orbit
$\smash{\{\widetilde{E}_1, \ldots, \widetilde{E}_6, \widetilde{E}'_1, \ldots, \widetilde{E}'_6\}}$
of size twelve is obtained from the first by applying the Geiser
involution~$g$.

Let~us now consider the auxiliary set
$\smash{\{E_1, \ldots, E_6, \widetilde{E}'_1, \ldots, \widetilde{E}'_6\}}$
of exceptional curves. Then,~according to formula~(\ref{projf}),

\begin{iii}
\item
$E_1,\ldots,E_6$
are mutually skew,
\item
$\smash{\widetilde{E}'_1, \ldots, \widetilde{E}'_6}$
are mutually skew, and
\item
$\smash{E_i \!\cdot\! \widetilde{E}'_j = 0}$,
for\/
$i \neq j$,
$1 \leq i,j \leq 6$,
and\/
$\smash{E_i \!\cdot\! \widetilde{E}'_i = 1}$
for~$i = 1, \ldots, 6$.
\end{iii}

\noindent
Thus,~Lemma~\ref{crit_hex} shows that
$$\smash{\{(\pi(E_1), \pi(\widetilde{E}'_1)), \ldots, (\pi(E_6), \pi(\widetilde{E}'_6))\} = \{(\pi(E_1), \pi(E'_1)), \ldots, (\pi(E_6), \pi(E'_6))\}}$$
is a Steiner~hexad. In~particular, the image
$\overline\iota(U_\ds)$
stabilises a Steiner hexad and is hence contained
in~$U_{63}$,
as claimed.
}
\eop
\end{theo}

\begin{rem}
The group
$W(E_6)$
has 350 conjugacy classes of~subgroups. Among these, 102 stabilise a double-six, i.e.\ are contained
in~$U_\ds$,
so that Theorem~\ref{dsix}~applies. For~each such conjugacy class, we constructed an example of a of cubic surface
over~$\bbQ$
using a different method~\cite{EJ10}.
\end{rem}

\section{Another application: Cubic surfaces with a rational line}
\label{sieben}

The maximal subgroup
$U_\l \subset G_\maxi$
of index~27 is just the stabiliser of a~line. We~have the following application to cubic surfaces with a rational line, which seems to be a slight refinement of the results given in~\cite[Section~6]{KST}.

\begin{theo}
Let an infinite
field\/~$k$
of characteristic
not\/~$2$,
a normal and separable extension
field\/~$l$,
and an injective group~homomorphism
$$i\colon \Gal(l/k) \hookrightarrow U_\l \subset G_\maxi \cong W(E_6) \subset S_{27}$$
be~given. Then~there exists a nonsingular cubic
surface\/~$S$
over\/~$k$
such~that

\begin{iii}
\item
the 27 lines
on\/~$S$
are defined
over\/~$l$
and each\/
$\sigma \in \Gal(l/k)$
permutes them as described by\/
$i(\sigma) \in U_\ds \subset S_{27}$.
\item
$S$
is\/
\mbox{$k$-unirational}.
\end{iii}\smallskip

\noindent
{\bf Proof.}
{\em
As~above, it suffices to show that
$\overline\iota(U_\l)$
is contained
in~$U_{63}$.
For~this, once again, we consider a model cubic surface over an algebraically closed~field. The~group
$U_\l$
operates on the 27~lines with an orbit type
$[1, 10, 16]$,
the orbit of size ten consisting of the lines that intersect the invariant
line~$L$.
Hence,~$\iota(U_\l)$
operates on the 56 exceptional curves of the degree two del Pezzo surface, obtained by blowing up one point, with orbit type
$[1,1,1,1,10,10,16,16]$.

It~is well known~\cite[Lemma~IV.15]{Be} that the exceptional curves in any of the size ten orbits may be written in the form
$E_1, \ldots, E_5, E'_1, \ldots, E'_5$,
where

\begin{iii}
\item
$E_1,\ldots,E_5$
are mutually skew,
\item
$E'_1,\ldots,E'_5$
are mutually skew, and
\item
$E_i \!\cdot\! E'_j = 0$,
for\/
$i \neq j$,
$1 \leq i,j \leq 5$,
and\/
$E_i \!\cdot\! E'_i = 1$
for~$i = 1, \ldots, 5$.
\end{iii}

\noindent
In~addition, the
image~$E_0$
of the invariant line fulfils
$E_0 \!\cdot\! E_i = E_0 \!\cdot\! E'_i = 1$,
for~$1, \ldots, 5$.
Finally,~the inverse
image~$E$
of the blow-up point is skew to all these~curves.

Formula~(\ref{projf}) now shows that
$E_1, \ldots, E_5, E, E'_1, \ldots, E'_5, \widetilde{E}_0$
fulfil, in this order, the assumptions of Lemma~\ref{crit_hex}. In~particular, the image
$\overline\iota(U_\l)$
stabilises a Steiner hexad and is hence contained
in~$U_{63}$,
as~required.
}
\eop
\end{theo}

\begin{coro}
Let an infinite
field\/~$k$
of characteristic
not\/~$2$,
a normal and separable extension
field\/~$l$,
and an injective group~homomorphism
$$i\colon \Gal(l/k) \hookrightarrow g_\maxi \cong W(D_5) \subset S_{16}$$
be~given. Then~there exists a del Pezzo
surface\/~$D$
of degree four
over\/~$k$
such~that

\begin{iii}
\item
the 16 exceptional curves
on\/~$D$
are defined
over\/~$l$
and each\/
$\sigma \in \Gal(l/k)$
permutes them as described by\/
$i(\sigma) \in g_\maxi \subset S_{16}$.
\item
$D$
is\/
\mbox{$k$-unirational}.
\eop
\end{iii}
\end{coro}

\setlength\parindent{0mm}

\frenchspacing
\begin{thebibliography}{EGA\,II}
\bibitem[Be]{Be}
Beauville, A.: Complex algebraic surfaces, LMS Lecture Note Series~68, {\em Cambridge University Press,} Cambridge~1983
\bibitem[BCP]{BCP}
Bosma, W., Cannon, J., and Playoust, C.: The Magma algebra system~I. The user language, {\em J.\ Symbolic Comput.} {\bf 24}\br(1997)\brr235--265
\bibitem[CHM]{CHM}
Conway, J.\,H., Hulpke, A., and McKay, J.: On transitive permutation groups, {\em LMS J.\ Comput.\ Math.} {\bf 1}\br(1998)\brr1--8
\bibitem[Co]{Co}
Corn, P.: The Brauer-Manin obstruction on del Pezzo surfaces of degree 2, {\em
Proc.\ Lond.\ Math.\ Soc.} {\bf 95}\br(2007)\brr735--777
\bibitem[De]{De}
Deligne, P.: Le groupe fondamental de la droite projective moins trois points, in: Galois groups over
$\bbQ$,
Berkeley~1987, Math.\ Sci.\ Res.\ Inst.\ Publ.~16, {\em Springer,} New York~1989
\bibitem[DM]{DM}
Dixon, J.\,D.\ and Mortimer, B.: Permutation groups, Graduate Texts in Mathematics~163, {\em Springer,} New York~1996
\bibitem[Do]{Do}
Dolgachev, I.\,V.: Classical Algebraic Geometry: a modern view, {\em Cambridge University press,} Cambridge~2012
\bibitem[DO]{DO}
Dolgachev, I.\,V.\ and Ortland, D.: Point sets in projective spaces and theta functions, {\em Ast\'erisque\/} {\bf 165}\br(1988)
\bibitem[EGA\,II]{EGAII}
Grothendieck, A.\ and Dieudonn\'e, J.: \smash{\'E}tude globale \'el\'ementaire de quelques classes de morphismes (EGA\,II), {\em Publ.\ Math.\ IHES\/} {\bf 8}\br(1961)
\bibitem[El]{El}
Elsenhans, A.-S.: Explicit computations of invariants of plane quartic curves, {\em Journal of Symbolic Computation\/} {\bf 68}\br(2015)\brr109--115
\bibitem[EJ10]{EJ10}
Elsenhans, A.-S. and Jahnel, J.: Cubic surfaces with a Galois invariant double-six, {\em Central European Journal of Mathematics\/} {\bf 8}\br(2010)\brr646--661
\bibitem[EJ15]{EJ15}
Elsenhans, A.-S. and Jahnel, J.: Moduli spaces and the inverse Galois problem for cubic surfaces, {\em Trans.\ AMS\/} {\bf 367}\br(2015)\brr7837--7861 
\bibitem[Er]{Er}
Ern\'e, R.: Construction of a del Pezzo surface with maximal Galois action on its Picard group, {\em J.\ Pure Appl.\ Algebra\/} {\bf 97}\br(1994)\brr15--27
\bibitem[Ha]{Ha}
Hartshorne, R.: Algebraic Geometry, Graduate Texts in Math.~52, {\em Springer,} New York, Heidelberg, Berlin~1977
\bibitem[KM]{KM}
Kl\"uners, J.\ and Malle, G.: Explicit Galois realization of transitive groups of degree up to 15, {\em J.\ Symbolic Comput.} {\bf 30}\br(2000)\brr675--716
\bibitem[Ko96]{Ko96}
Koll\'ar, J.: Rational curves on algebraic varieties, Ergebnisse der Mathematik und ihrer Grenzgebiete~32, {\em Springer,} Berlin~1996
\bibitem[Ko02]{Ko02}
Koll\'ar, J.: Unirationality of cubic hypersurfaces, {\em J.\ Inst.\ Math.\ Jussieu\/} {\bf 1}\br(2002)\brr467--476
\bibitem[KST]{KST}
Kunyavskij, \smash{B.\,\`E.,} Skorobogatov, A.\,N., and Tsfasman, M.\,A.: Del Pezzo surfaces of degree four, {\em M\'em.\ Soc.\ Math.\ France\/} {\bf 37}\br(1989)\brr1--113
\bibitem[Ma]{Ma}
Manin, Yu.\,I.: Cubic forms, algebra, geometry, arithmetic,
{\em North-Holland Publishing Co.} and {\em American Elsevier Publishing Co.,}
Amsterdam, London, and New York~1974
\bibitem[Ne]{Ne}
Neukirch, J.: Algebraic number theory, Grundlehren der Mathematischen Wis\-sen\-schaf\-ten~322, {\em Springer,} Berlin~1999
\bibitem[NSW]{NSW}
Neukirch, J., Schmidt, A., and Wingberg, K.: Cohomology of number fields, Grundlehren der Mathematischen Wissenschaften~323, {\em Springer,} Berlin~2000
\bibitem[Pl]{Pl}
Pl\"ucker, J.: Solution d'une question fondamentale concernant la th\'eorie g\'en\'erale des courbes, {\em J.~f\"ur die reine und angew.\ Math.} {\bf 12}\br(1834)\brr105--108
\bibitem[Ro]{Ro}
Royle, G.\,F.: The transitive groups of degree twelve, {\em J.\ Symbolic Comput.} {\bf 4}\br(1987)\brr255--268
\bibitem[SGA1]{SGA1}
Grothendieck, A.: Rev\^etements \'etales et groupe fondamental (SGA\,1), Lecture Notes Math.~224, {\em Springer,} Berlin~1971
\bibitem[Sha]{Sha}
Shafarevich, I.\,R.: Construction of fields of algebraic numbers with given solvable Galois group (Russian), {\em Izv.\ Akad.\ Nauk SSSR\/} {\bf 18}\br(1954)\brr525--578
\bibitem[Shi]{Shi}
Shioda, T.: Plane quartics and Mordell-Weil lattices of type
$E_7$,
{\em Comment.\ Math.\ Univ.\ St.\ Pauli\/} {\bf 42}\br(1993)\brr61--79
\bibitem[Sk13]{Sk13}
Skorobogatov, A.\,N.: Arithmetic Geometry: Rational Points, {\em Preprint\/} available at {\tt https:\discretionary{}{}{}//wwwf.imperial.ac.uk/\~{}anskor/arith\_geom\_files/TCCnotes.pdf}
\bibitem[Sk17]{Sk17}
Skorobogatov, A.\,N.: Cohomology and the Brauer group of double covers, in: Brauer groups and obstruction problems, Progr. Math.~320, {\em Birkh\"auser/Springer,} Cham~2017, 231--247
\bibitem[SV]{SV}
St\"ohr, K.-O.\ and Voloch, J.\,F.: A formula for the Cartier operator on plane algebraic curves, {\em J.~f\"ur die reine und angew.\ Math.} {\bf 377}\br(1987)\brr49--64
\bibitem[Va]{Va}
V\'arilly-Alvarado, A.: Arithmetic of del Pezzo surfaces, in: Birational geometry, rational curves, and arithmetic, {\em Springer,} New York~2013, 293--319
\bibitem[Za]{Za}
Zarhin, Yu.\,G: Del Pezzo surfaces of degree 1 and Jacobians, {\em Math.\ Ann.} {\bf 340}\br(2008)\brr407--435
\end{thebibliography}
\end{document}